\documentclass[a4paper,10pt]{article}
\usepackage{amsmath}
\usepackage{amssymb,amsfonts,amsthm}
\usepackage[utf8]{inputenc} 
\usepackage[english]{babel}
\usepackage{lpic}
\usepackage{anysize}
\usepackage{enumerate}
\usepackage{tikz}
\usepackage{subfig}
\usepackage[all]{xy}
\usepackage{mathrsfs}
\usepackage{verbatim}
\usepackage{cite}
\definecolor{green}{rgb}{0,0.8,0.5}

\usepackage[pdftex]{hyperref}
\hypersetup{colorlinks=true,linkcolor=blue,citecolor=red}
\usepackage[width=0.8\textwidth]{caption}

\renewenvironment{abstract}{\small\quotation\noindent
 {\bfseries \abstractname .}}{\endquotation \par}

\newenvironment{trivlistparm}[2]{\trivlist
\item[{#1 #2}]}{\endtrivlist}

\newenvironment{prooftext}[1]{\trivlistparm{\bfseries}{#1}}{\Qed\endtrivlistparm}
\newenvironment{prova}{\trivlistparm{\bfseries}{Proof.}}{\Qed\endtrivlistparm}

\catcode`\@=11

\def\resetthefootnote{\renewcommand{\thefootnote}{\@arabic\c@footnote} }
\def\@principiremex#1{\trivlist
 \item[\hskip \labelsep{\bfseries #1\ \thethm.}]\ignorespaces}
\def\opar@principiremex#1[#2]{\trivlist
 \item[\hskip \labelsep{\bfseries #1\ \thethm\ (#2).}]\ignorespaces}

\newcommand{\newTHEOremrom}[2]{\newenvironment{#1}{\refstepcounter{thm}\@ifnextchar[{\opar@principiremex{#2}}
{\@principiremex{#2}}}{\qedB\endtrivlist}} \catcode`\@=12
\DeclareMathSymbol{\square}{\mathord}{AMSa}{"03}
\newcommand{\qedB}{\nopagebreak\hspace*{\fill}$\square$\par}
\newcommand{\Qed}{\nopagebreak\hspace*{\fill}{\vrule width6pt height6pt depth0pt}\par}

\newTHEOremrom{defi} {Definition}
\newTHEOremrom {rem} {Remark}
\newTHEOremrom {ex} {Example}

\renewcommand{\geq}{\geqslant}
\renewcommand{\epsilon}{\varepsilon}
\renewcommand{\leq}{\leqslant}
\newcommand{\F}{\mathscr{F}}
\newcommand{\LL}{\mathscr{L}}
\newcommand{\CC}{\mathscr{B}}

\newcommand{\DD}{\mathscr{D}}

\newtheorem{thm}{Theorem}[section]
\newtheorem{thmx}{Theorem}
\newtheorem{prop}[thm]{Proposition}

\newtheorem{lema}[thm]{Lemma}

\newcommand{\abs}[1]{\left|#1\right|}
\newcommand{\R}{\mathbb{R}}

\newcommand{\Z}{\mathbb{Z}}
\newcommand{\N}{\mathbb{N}}

\newcommand{\PA}{\mathscr{P}}
\newcommand{\PI}{\mathcal{I}}

\newcommand{\RP}{\mathbb{RP}}

\newcommand{\out}{\Pi}
\newcommand{\fun}{\mathscr G}

\newcommand{\sist}[2]{
  \left\{\!
   \begin{array}{l}
    \dot x=#1 \\[2pt] \dot y=#2
   \end{array}
  \right.
}

\title{\textbf{Asymptotic development of an integral operator and boundedness of the criticality of potential centers}
}

\author{David Rojas\\[10pt]
{\small \textsl{Departament d’Inform\`atica, Matem\`atica Aplicada i Estad\'istica,}}\\
\vspace{-2pt}
{\small \textsl{Universitat de Girona, 17003 Girona, Spain.}}\\[5pt]
}

\date{}

\begin{document}

\maketitle

\begin{abstract}
We study the asymptotic development at infinity of an integral operator. We use this development to give sufficient conditions to upper bound the number of critical periodic orbits that bifurcate from the outer boundary of the period function of planar potential centers. We apply the main results to two different families: the power-like potential family $\ddot x=x^p-x^q$, $p,q\in\R$, $p>q$; and the family of dehomogenized Loud's centers.
\end{abstract}

\section{Introduction}\label{sec:intro}

Consider a continuous family of planar potential systems
\begin{equation}\label{centro}
\dot x = -y,\ \dot y = V_{\mu}'(x),
\end{equation}
where $\mu\in\Lambda$ is a parameter, $\Lambda$ is an open subset of $\R^d$, $d\geq 1$, and $V_{\mu}$ is an analytic function defined in an open interval $I_{\mu}\subset\R$ containing $x=0$. In the case when $V_{\mu}(0)=V_{\mu}'(0)=0$ and $V_{\mu}''(0)>0$ equation~\eqref{centro} has a non-degenerate center at the origin for each value of the parameter and so the point $(0,0)$ has a punctured neighbourhood that is entirely foliated by periodic orbits surrounding it. The largest neighbourhood with this property is the period annulus of the center and we shall denote it by $\PA_{\mu}$. If we consider the embedding of $\PA_{\mu}$ into $\RP^2$, its boundary, namely $\partial\PA_{\mu}$, is divided in two connected components: the origin itself, which is called the inner boundary of the period annulus, and the outer boundary of the period annulus defined by $\Pi_{\mu}\!:=\partial \PA_{\mu}\setminus\{(0,0)\}$. When the center is a potential oscillator the natural parametrization of the closed orbits inside the period annulus is given by the energy level of the Hamiltonian $H(x,y;\mu)=\tfrac{1}{2}y^2+V_{\mu}(x)$. Since $V_{\mu}(0)=0$ by convention, we have that $H(\PA_{\mu})=(0,h_0(\mu))$, where $h_0(\mu)\in \R^+\cup\{+\infty\}$ denotes the energy level of the outer boundary $\Pi_{\mu}$.

The object under study in this paper is the period function of the center. The minimal period $T_{\mu}(h)$ of the periodic orbit $\gamma_{h,\mu}$ inside the energy level $\{H(x,y;\mu)=h\}$ can be written as the Abelian integral
\[
T_{\mu}(h)=\int_{\gamma_{h,\mu}} \frac{dx}{y}.
\]
This function is analytic on $(0,h_0(\mu))$ for each value of the parameter and it can be extended analytically to $h=0$ due to the non-degeneracy of the center. The derivative $T_{\mu}'(h)$ can also be written as an Abelian integral and its zeros correspond to critical periodic orbits of the system. This paper is concerned with the bifurcation of such critical periodic orbits from the outer boundary $\Pi_{\mu}$. That is, for a fixed $\mu_0\in\Lambda$, we aim to control the number of critical periodic orbits of system~\eqref{centro} that may emerge or disappear from $\Pi_{\mu_0}$ as we move slightly the parameter $\mu\approx\mu_0$. This number is called the criticality of the outer boundary.

\begin{defi}\label{defi:criticality}
Consider a continuous family $\{X_{\mu}\}_{\mu\in\Lambda}$ of planar analytic vector fields with a center and fix some $\mu_0\in\Lambda$. Suppose that the outer boundary of the period annulus varies continuously at $\mu_0\in\Lambda$, meaning that for any $\epsilon>0$ there exists $\delta>0$ such that $d_H(\Pi_{\mu},\Pi_{\mu_0})\leq \epsilon$ for all $\mu\in\Lambda$ with $\|\mu-\mu_0\|\leq \delta$. Then, setting
\[
N(\delta,\epsilon)\!:=\sup\{\# \text{critical periodic orbits }\gamma\text{ of }X_{\mu}\text{ in }\PA_{\mu}\text{ with }d_H(\gamma,\Pi_{\mu_0})\leq \epsilon \text{ and }\|\mu-\mu_0\|\leq \delta  \},
\]
the criticality of $(\Pi_{\mu_0},X_{\mu_0})$ with respect to the deformation $X_{\mu}$ is $\mathrm{Crit}\bigl((\Pi_{\mu_0},X_{\mu_0}),X_{\mu}\bigr)\!:=\inf_{\delta,\epsilon} N(\delta,\epsilon)$.
\end{defi}

In the previous definition $d_H$ stands for the Hausdorff distance between compact sets of $\RP^2$. Notice that according with this definition the criticality may be infinite but, in the case it is not, it gives the maximal number of critical periodic orbits of $X_{\mu}$ tending to the outer boundary $\Pi_{\mu_0}$ in the Hausdorff sense as the parameter $\mu$ approaches $\mu_0$. The requirement of the continuity of $\PA_{\mu}$ with respect to the parameters of the system ensures that the possible changes of $\PA_{\mu}$ do not occur abruptly. We refer to~\cite{ManVil2006} for details illustrating the necessity of this extra assumption. 

\begin{defi}
A parameter $\mu_0\in\Lambda$ is called a local regular value of the period function at the outer boundary of the period annulus if $\mathrm{Crit}\bigl((\Pi_{\mu_0},X_{\mu_0}),X_{\mu}\bigr)=0$. Otherwise the parameter is called a local bifurcation value at the outer boundary.
\end{defi}

The present paper is a contribution that follows the spirit of the series of works~\cite{ManRojVil2016,ManRojVil2016a,Rojas2018}. In these papers, we develop analytical tools which allow to give an upper bound of the criticality at the outer boundary of the period annulus of families of planar potential systems~\eqref{centro}. The key idea is to find a collection of functions $\phi_{\mu}^i(h)$, $i=1,2,\dots,n$, verifying that there exist $\delta,\epsilon>0$ such that $(\phi_{\mu}^1,\phi_{\mu}^2,\dots,\phi_{\mu}^n, T_{\mu}')$ form an Extended Complete Chebyshev system (ECT-system for short, see Definition~\ref{defi_ECT}) on the interval $(h_0(\mu)-\epsilon,h_0(\mu))$ and $\|\mu-\mu_0\|\leq \delta$. This fact implies that $T_{\mu}'(h)$ has at most $n$ zeros in $(h_0(\mu)-\epsilon,h_0(\mu))$, counting multiplicities, uniformly on the parameters $\mu\approx\mu_0$. 
In particular, $\mathrm{Crit}\bigl((\Pi_{\mu_0},X_{\mu_0}),X_{\mu}\bigr)\leq n$. According with Lemma~\ref{lema:ECT-Wronskia}, to give an upper bound of the criticality is reduced to guarantee that the Wronskian (see Definition~\ref{wronskian}) $W[\phi_{\mu}^1,\phi_{\mu}^2,\dots,\phi_{\mu}^n, T_{\mu}'](h)$ does not vanish for all $(h,\mu)\approx(h_0(\mu_0),\mu_0)$. The tools of the previous works, and also the ones we present here, allow to tackle this problem in the following two situations: either $h_0(\mu)=+\infty$ or $h_0(\mu)<+\infty$ for all $\mu\approx\mu_0$. That is, the case in which there exist $\mu_1$ and $\mu_2$ in any neighbourhood of $\mu_0$ with $h_0(\mu_1)=+\infty$ and $h_0(\mu_2)<+\infty$ is not considered.

Roughly speaking, the results in the previous papers relate the first term in the asymptotic development of the potential $V_{\mu}$ at the endpoints of $I_{\mu}$ with the first term of the asymptotic development of $W[\phi_{\mu}^1,\phi_{\mu}^2,\dots,\phi_{\mu}^n, T_{\mu}'](h)$ at $h=h_0(\mu)$. This allows to control the sign of the Wronskian under consideration for $(h,\mu)\approx(h_0(\mu_0),\mu_0)$; that is, uniformly on the parameters $\mu\approx\mu_0$. However, there are some situations where these first terms of $V_{\mu}$ are not enough to compute the first term of the Wronskian at $h=h_0(\mu)$ and so more terms in the asymptotic development must be employed. Theorem~\ref{thm:criticalitat-infinite} and~\ref{thm:criticalitat-finite} in Section~\ref{sec:dynamic} aim to generalize the results in~\cite{ManRojVil2016,ManRojVil2016a,Rojas2018} in this direction. To accomplish the desired results, we will employ a generalization of~\cite[Proposition 2.16]{ManRojVil2016a} and~\cite[Theorem D]{Rojas2018} (see Theorem~\ref{thm:tecnic} in Section~\ref{sec:technical}.)

\begin{figure}
\centering
\subfloat[Previous bifurcation diagram\label{fig:familia-abans}]
{
 \includegraphics[scale=1]{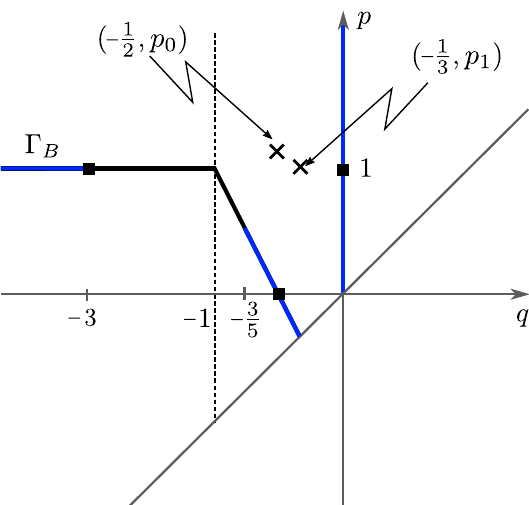}
}
\hspace{1cm}
\subfloat[New bifurcation diagram\label{fig:familia-nou}]
{
 \includegraphics[scale=1]{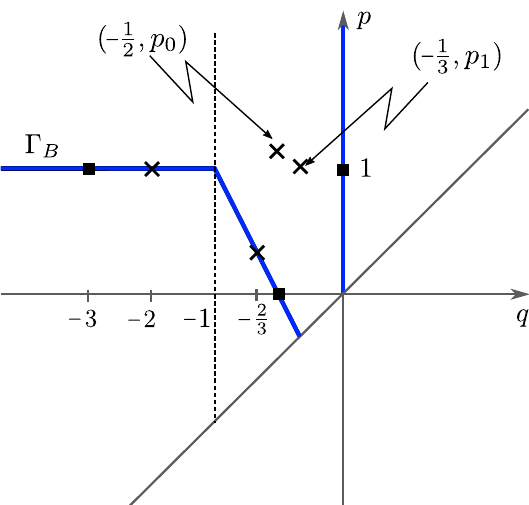}
}
\caption{On the left, bifurcation diagram of the period function of the family~\eqref{familia-potencial} at the outer boundary of the period annulus according with \cite{ManRojVil2016,ManRojVil2016a,Rojas2018}. On the right, improvement of the bifurcation diagram according with Theorem~\ref{thm:familia}. In both figures, $\Gamma_B$ stands for the union of the bold lines. In black the parameters with criticality at least one. In blue the parameters with criticality exactly one. Black squares are parameters that correspond to isochronous centers. Crosses are parameters where techniques do not apply. The $q$-axis do not preserve the scale for the sake of space. 
}
\end{figure}

As an illustration of these generalizations we recover the study of two different families of planar centers. The first application is on the two-parametric family of potential differential system given by
\begin{equation}\label{familia-potencial}
\sist{-y,}{(x+1)^p-(x+1)^q,}
\end{equation}
which has a non-degenerate center at the origin for all $\mu\!:=(q,p)$ varying in $\Lambda\!:=\{(q,p)\in\R^2: p>q\}$. As far as we know, the period function of this center was originally studied by Miyamoto and Yagasaki~\cite{MY}, giving a monotonicity result for $q=1$ and $p\in\N$. This result was improved later by Yagasaki~\cite{Y} showing that the period function of~\eqref{familia-potencial} is monotonous for $q=1$ and any $p>1$ real. Motivated by those results, we considered the whole family~\eqref{familia-potencial} with $p>q$ and performed an exhaustive study of the period function in~\cite{ManRojVil2017}. Concerning the criticality at the outer boundary, the family~\eqref{familia-potencial} became our testing ground for the techniques mentioned before. In these works, see Figure~\ref{fig:familia-abans}, we proved that parameters $\mu_0\in\Lambda\setminus\bigl( \Gamma_B\cup\{q+1=0\}\cup\{(-\tfrac{1}{2},p_1)\}\cup\{(-\tfrac{1}{3},p_2)\}\bigr)$, with $p_1\approx 1.20175$ and $p_2\approx 1.15685$, are local regular values of the period function at the outer boundary. In addition, $\mathrm{Crit}\bigl((\Pi_{\mu_0},X_{\mu_0}),X_{\mu}\bigr)\geq 1$ if $\mu_0\in\Gamma_B$ and it is exactly one for parameters $\mu_0=(q_0,p_0)$ satisfying either $q_0=0$ and $p_0\in(0,+\infty)\setminus\{1\}$, $p_0=1$ and $q_0<-3$, or $p_0+2q_0+1=0$ and $q_0\in(-\tfrac{3}{5},-\tfrac{1}{3})\setminus\{-\tfrac{1}{2}\}$. Using the tools in the present paper, the bifurcation diagram in Figure~\ref{fig:familia-abans} is improved by the following result. (See Figure~\ref{fig:familia-nou}.)

\begin{thmx}\label{thm:familia}
Let $\{X_{\mu}\}_{\mu\in\Lambda}$ be the family of analytic potential systems~\eqref{familia-potencial} and consider the period function of the center at the origin. If $\mu_0=(q_0,p_0)$ with either $p_0=1$ and $q_0\in(-3,-1)\setminus\{-2\}$, or $p_0+2q_0+1=0$ and $q_0\in(-1,-\tfrac{1}{2})\setminus\{-\tfrac{2}{3}\}$ then $\mathrm{Crit}\bigl((\Pi_{\mu_0},X_{\mu_0}),X_{\mu}\bigr)=1$.
\end{thmx}

The proof of this result is presented in Section~\ref{sec:lambda2} (see Proposition~\ref{prop:fam1}) and Section~\ref{sec:lambda1} (see Proposition~\ref{prop:resultat-familia}). The result finishes the bifurcation diagram of the period function at the outer boundary of the family~\eqref{familia-potencial} except for the line $\{q+1=0\}$ and four points in the parameter space. (See Figure~\ref{fig:familia-nou}.) We point out that the line corresponds to parameters such that the energy at the outer boundary $h_0(\mu)$ changes from infinite ($q<-1$) to finite ($q>-1$). The points correspond to parameters that do not satisfy the technical hypothesis to apply the analytic tools. The bifurcation diagram in Figure~\ref{fig:familia-nou} agrees with the global bifurcation diagram conjectured in~\cite{ManRojVil2017}.

\begin{figure}
\centering
\includegraphics[scale=.8]{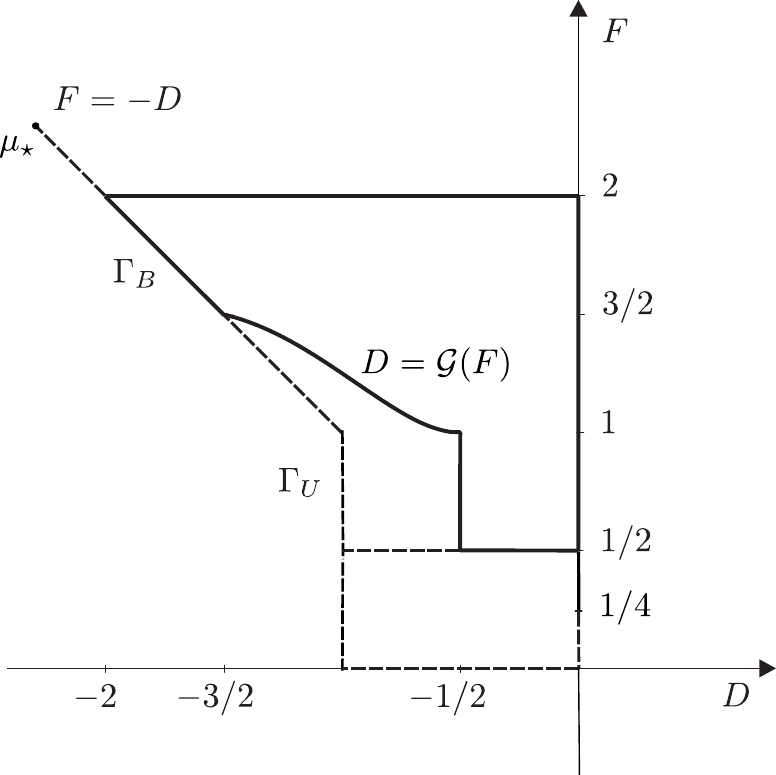}
\caption{\label{fig:loud}Bifurcation diagram of the period function at the polycycle of system~\eqref{loud}, where $\mu_{\star}=(-F_{\star},F_{\star})$ with $F_{\star}\approx 2.34$. The union of the bold curves correspond to the set of bifurcation parameters at the outer boundary. The union of dotted straight lines correspond to the set of unspecified parameters. The complementary of those two sets correspond to regular parameters.}
\end{figure}

The second application is on the family of quadratic polynomial planar centers. The literature classify quadratic centers in four families: Hamiltonian, reversible $Q_3^R$, codimension four $Q_4$, and generalized Lotka-Volterra $Q_3^{LV}$. Chicone~\cite{Chicone} conjectured that reversible centers have at most two critical periodic orbits whereas the centers of the other three families have monotonic period function. Regarding quadratic reversible centers, by an affine transformation and a constant rescaling of time, they can be brought to the Loud normal form
\[
\sist{-y+Bxy,}{x+Dx^2+Fy^2.}
\]
In~\cite{GGV} the authors show that if $B=0$ the period of the center at the origin is globally monotone. When $B\neq0$ one can reduce the system, by means of a rescaling, to $B=1$. That is,
\begin{equation}\label{loud}
\sist{-y+xy,}{x+Dx^2+Fy^2.}
\end{equation}
This family is known as dehomogenized Loud's centers and it has a center at the origin for all parameters $\mu\!:=(D,F)\in\R^2$. The bifurcation of critical periodic orbits from the outer boundary of the period annulus of system~\eqref{loud} has been extensively studied in the recent years. (See Figure~\ref{fig:loud}.) We refer to the series of papers \cite{ManVil2006,MarMarSaaVil2015,MarMarVil2006,MarVil2006,Villadelprat2007,RojVil2018} and references therein. 

Our contribution to the bifurcation diagram at the outer boundary of the dehomogenized Loud's centers is to show that almost all parameters in the bifurcation curve $D=\mathcal G(F)$ have criticality exactly one. Up to now, this have been proved for parameters in that bifurcation curve with $F\in(\tfrac{3}{2},\tfrac{4}{3})$. In order to state the result properly, let us consider the parameter space 
\[
\Lambda\!:=\{(D,F)\in\R^2 : 1<F<\tfrac{3}{2}, D<-\tfrac{1}{2}, D+F>0\}.
\]
Moreover, let ${}_2F_1(a,b;c;z)$ be the Hypergeometric function (see~\cite[Section~15]{AS}), $\alpha\!:=(p_2-1)/(p_2-p_1)$ with $p_1$ and $p_2$ defined in~\eqref{def:p1p2}, and let $c(\mu)$ be the function in the statement of Lemma~\ref{lema:quant-S-dreta} using the expression of $V_{\mu}$ in~\eqref{V-loud}.

\begin{thmx}\label{thm:loud}
Let $\{X_{\mu}\}_{\mu\in\Lambda}$ be the family of analytic potential systems~\eqref{loud} and consider the period function of the center at the origin. Let $\mu_0=(D_0,F_0)\in\Lambda$ satisfying $D_0=\mathcal G(F_0)$. Then $\mathrm{Crit}\bigl((\Pi_{\mu_0},X_{\mu_0}),X_{\mu}\bigr)=1$ in the following situations:
\begin{enumerate}[$(a)$]
\item If $F_0\in(\tfrac{6}{5},\tfrac{4}{3})$.
\item If $F_0\in(\tfrac{9}{8},\tfrac{6}{5})$ and $c(\mu_0)\neq 0$.
\item If $F_0\in(1,\tfrac{9}{8}]$, $c(\mu_0)\neq 0$ and 
\begin{equation}\label{eq:condition}
\frac{{}_2F_1\bigl(-\frac{3}{2},\frac{5}{2};\frac{7}{2}-4F_0;\alpha(\mu_0)\bigr)}{\Gamma\bigl(\tfrac{7}{2}-4F_0\bigr)}\neq 0.
\end{equation}
\end{enumerate}
\end{thmx}

The proof of the result is given in Sections~\ref{sec:proof_ThmB} and~\ref{sec:comments}. In that last Section, also a numerical manifestation that condition~\eqref{eq:condition} seems to be fulfilled for $F\in(1,\tfrac{9}{8}]$ is given. The condition $c(\mu_0)\neq 0$ is a technical requirement for the techniques involved in the proof and it is conjectured to be non necessary for the criticality to be also exactly one for those parameters satisfying $c(\mu_0)=0$. Regarding the equality $c(\mu_0)=0$ we also show numerically that the equation has a unique solution $(D_0,F_0)\approx(-0.56996,1.00781)$ in Section~\ref{sec:comments}. The parameter $D_0=\mathcal G(\tfrac{4}{3})$ is conjectured to have criticality exactly two at the outer boundary of the period annulus.

In the forthcoming paper~\cite{MarVilpreprint} the authors also obtain~\eqref{eq:condition} as a sufficient condition for bifurcation parameters in $D=\mathcal G(F)$ to have criticality exactly one using a completely different approach. In that work the authors study the asymptotic development of the Dulac time function near hyperbolic saddle singularities of meromorphic planar centers. They also use the dehomogenized Loud's centers as testing ground and reach the same result in Theorem~\ref{thm:loud} without the technical restriction $c(\mu_0)\neq 0$ and allowing $F_0=\tfrac{6}{5}$.

The rest of the paper is organized as follows. In Section~\ref{sec:technical} we study the asymptotic behaviour of an integral operator that will be useful for the proof of the dynamical results in the paper. In Section~\ref{sec:dynamic} we use these techniques to obtain two results that gives sufficient conditions in order to bound the criticality at the outer boundary of planar potential centers. Finally Section~\ref{sec:applications} is dedicated to the application of such sufficient conditions and so we prove Theorems~\ref{thm:familia} and~\ref{thm:loud}. The work is complemented with an Appendix that contain the more technical proofs.

\section{Asymptotic behaviour of a certain integral operator}\label{sec:technical}
Let $c\in\R^+\cup\{+\infty\}$ and consider the integral operator
\[
\F:C[0,c) \rightarrow C[0,c)
\]
defined by
\[
\F[f](x)\!:= \int_0^{\frac{\pi}{2}} f(x\sin\theta)d\theta.
\]
Here, and in what follows, $C[0,c)$ stands for the set of continuous functions on $[0,c)$. In the collection of works \cite{ManRojVil2016,ManRojVil2016a,Rojas2018} the previous operator is studied because of its relation with the bifurcation of critical periodic orbits. Indeed, the derivative of the period function of system~\eqref{centro} satisfies the equality
\begin{equation}\label{FT}
\sqrt{2}h^2 T_{\mu}'(h^2)=\F[f_{\mu}](h),\ h\in(0,h_0(\mu))
\end{equation}
with $f_{\mu}(x)=x(g_{\mu}^{-1})''(x)-x(g_{\mu}^{-1})''(-x)$ and $g_{\mu}(x)\!:=\text{sgn}(x)\sqrt{V_{\mu}(x)}$. Roughly speaking, the main objective in these works is to give sufficient conditions to the first term of the asymptotic development at $x=c$ of the function $f_{\mu}$ in order that the first term of the asymptotic development at $x=c$ of $\F[f_{\mu}]$ is obtained, uniformly on the parameters. These conditions are formulated using the following notions.

\begin{defi}\label{def:cont_quantifiable}
Let $\{f_{\mu}\}_{\mu\in\Lambda}$ be a continuous family of continuous
functions on~$\bigl(a(\mu),b(\mu)\bigr)$, meaning that the map $(x,\mu)\longmapsto f_{\mu}(x)$ is continuous on $\bigl\{(x,\mu)\in\R\times\Lambda: x\in \bigl(a(\mu),b(\mu)\bigr)\bigr\}$. Assume that either $b:\Lambda \rightarrow\R$ is continuous or $b\equiv +\infty$ in $\Lambda$. Given $\mu_0\in\Lambda$ we say that $\{f_{\mu}\}_{\mu\in\Lambda}$ is
\emph{continuously quantifiable} in $\mu_0$ at $b(\mu)$ by
$\alpha(\mu)$ with limit $\ell(\mu)$ if there exists an open neighbourhood~$U$ of $\mu_0$ such that for all $\hat\mu\in U$,
\begin{enumerate}[$(i)$]
\item If $b(\mu_0)<+\infty$, then $\lim_{(x,\mu)\rightarrow (b(\hat\mu),\hat\mu)} f_{\mu}(x)(b(\mu)-x)^{\alpha(\mu)}=\ell(\hat\mu)$ and $\ell(\hat\mu)\neq 0$.
\item If $b(\mu_0)=+\infty$, then $\lim_{(x,\mu)\rightarrow (+\infty,\hat\mu)} x^{-\alpha(\mu)}f_{\mu}(x)=\ell(\hat\mu)$ and $\ell(\hat\mu)\neq 0$.
\end{enumerate}
For the sake of shortness, in the first case we write $f_{\mu}(x)\sim_{b(\mu)} \ell(\mu)(b(\mu)-x)^{-\alpha(\mu)}$ at $\mu_0$, and in the second case $f_{\mu}(x)\sim_{+\infty} \ell(\mu)x^{\alpha(\mu)}$ at $\mu_0$. We use the analogous definition for the left endpoint $a(\mu)$.
\end{defi}

We point out that the map $\alpha:U\rightarrow\R$ in the previous definition is continuous at $\mu=\mu_0$ (see \cite[Remark~2.6]{ManRojVil2016}).

From now on let us assume that $f_{\mu}\in C[0,+\infty)$. The purpose of this section is to deal with an specific situation that was not contemplated in the previous works. With this aim in view we first recover the main results in~\cite{ManRojVil2016,Rojas2018}.

\begin{defi}
The function defined for all $x>0$ and $\alpha\in\R$ by means of
\[
\omega(x,\alpha)\!:=\begin{cases}
\frac{x^{\alpha+1}-1}{\alpha+1} & \text{ if }\alpha\neq -1,\\
\log x & \text{ if }\alpha=-1,
\end{cases}
\]
is called the \emph{Roussarie-Ecalle} compensator. For the sake of brevity, we also define
\[
\fun(\alpha)\!:=\frac{\sqrt{\pi}}{2}\frac{\Gamma\left(\frac{1+\alpha}{2}\right)}{\Gamma\left(1+\frac{\alpha}{2}\right)}\ \text{ and }\ \Omega(x,\alpha)\!:=(\alpha+1)\fun(\alpha)\omega(x,\alpha),
\]
where $\Gamma$ is the Gamma function. Following the notation in Definition~\ref{def:cont_quantifiable}, we write $f_{\mu}(x)\sim_{+\infty} \ell(\mu)\Omega(x,\alpha(\mu))$ at $\mu_0$ if 
\[
\lim_{(x,\mu)\rightarrow(+\infty,\mu_0)}\frac{f_{\mu}(x)}{\Omega(x,\alpha(\mu))}=\ell(\mu_0)\neq 0.
\]
\end{defi}

\begin{defi}
Let $f\in C[0,+\infty)$. We call
\[
M_n[f]\!:=\int_0^{+\infty}x^{2n-2} f(x)dx
\]
the $n$-th momentum of $f$, whenever it is well defined. If $n=1$ we simply say that $M[f]\!:=M_1[f]$ is the momentum of $f$.
\end{defi}

The next result gathers Theorems~$2.13$ and~$2.17$ in~\cite{ManRojVil2016} together with Theorem~C in \cite{Rojas2018}.

\begin{thm}\label{thm:quantificar}
Let $\Lambda$ be an open subset of $\R^d$ and consider a continuous family $\{f_{\mu}\}_{\mu\in\Lambda}$ of continuous functions on $[0,+\infty)$. Suppose that $f_{\mu}(x)\sim_{+\infty} a(\mu)x^{\alpha(\mu)}$ at $\mu_0$. The following assertions hold:
\begin{enumerate}[$(a)$]
\item If $\alpha(\mu_0)>-1$ then $\F[f_{\mu}](x)\sim_{+\infty} a(\mu)\fun(\alpha(\mu))x^{\alpha(\mu)}$ at $\mu_0$.
\item If $\alpha(\mu_0)=-1$ then $\F[f_{\mu}](x)\sim_{+\infty} a(\mu)\Omega(x,\alpha(\mu))\frac{1}{x}$ at $\mu_0$.
\item If $\alpha(\mu_0)<-1$ let us take $m\in\N$ such that $\alpha(\mu_0)+2m\in[-1,1)$. In this case:
\begin{enumerate}[$(c1)$]
\item If $M_i[f_{\mu}]\equiv 0$ for $i=1,\dots,\ell-1$ and $M_{\ell}[f_{\mu_0}]\neq 0$ for some $1\leq \ell\leq m$, then 
\[
\F[f_{\mu}](x)\sim_{+\infty} M_{\ell}[f_{\mu}]x^{1-2\ell} \text{ at }\mu_0.
\]
\item If $M_i[f_{\mu}]\equiv 0$ for $i=1,\dots,m$ and $\alpha(\mu_0)+2m\notin\{-1,0\}$ then 
\[
\F[f_{\mu}](x)\sim_{+\infty} a(\mu)\prod_{i=1}^m\frac{\alpha(\mu)+2i}{\alpha(\mu)+2i-1}\fun(\alpha(\mu)+2m)x^{\alpha(\mu)} \text{ at }\mu_0.
\]
\item If $M_i[f_{\mu}]\equiv 0$ for $i=1,\dots,m$ and $\alpha(\mu_0)+2m=-1$ then 
\[
\F[f_{\mu}](x)\sim_{+\infty} a(\mu)\prod_{i=1}^m\frac{\alpha(\mu)+2i}{\alpha(\mu)+2i-1}\Omega(x,\alpha(\mu)+2m)x^{-2m-1} \text{ at }\mu_0.
\]
\end{enumerate}
\end{enumerate}
\end{thm}

In the previous result the first term in the asymptotic development at infinity of the function $\F[f_{\mu}](x)$ is provided as soon as $f_{\mu}(x)\sim_{+\infty} a(\mu)x^{\alpha(\mu)}$ at $\mu_0$ except in the case when $\alpha(\mu_0)=-2m$ for some $m\in\N$, $m\geq 1$, and $M_i[f]=0$ for $i=1,\dots,m$. In this special situation the hypothesis of $f_{\mu}$ to be quantifiable by $\alpha(\mu)=-2m$ is not enough to quantify $\F[f_{\mu}]$ at infinity. The following three functions exemplify this phenomena for $n=1$ even in the non-parametric scenario:
\begin{equation}\label{examples}
f(x)=\begin{cases}
x^{-2}+x^{-\frac{5}{2}} & x\geq 1,\\
\frac{22}{3}x-\frac{16}{3} & x\in[0,1),
\end{cases} \ 
g(x)=\begin{cases}
x^{-2}+x^{-5} & x\geq 1,\\
\frac{13}{2}x-\frac{9}{2} & x\in[0,1),
\end{cases} \ 
h(x)=\begin{cases}
x^{-2}+x^{-3} & x\geq 1,\\
7x-5 & x\in[0,1).
\end{cases}
\end{equation}
All these functions are quantifiable by $\alpha=-2$ at infinity and it is a computation to show that their momenta vanish. One can verify that $\F[f]$ and $\F[g]$ are quantifiable at infinity by $-\frac{5}{2}$ and $-3$ respectively, and that
\[
\lim_{x\rightarrow+\infty}\frac{x^3}{\log x}\F[h](x)=\frac{1}{2}.
\]
These three examples subscribe the idea that more information on the asymptotic development of the function $f_{\mu}$ is needed to quantify $\F[f_{\mu}]$ in this situation. To address this problem we generalize Definition~\ref{def:cont_quantifiable} as follows.

\begin{defi}\label{defi:asympt}
Let $\{f_{\mu}\}_{\mu\in\Lambda}$ be a continuous family of continuous functions on $(a(\mu),b(\mu))$. Assume that either $b:\Lambda\rightarrow\R$ is continuous or $b\equiv +\infty$ in $\Lambda$. Given $\mu_0\in\Lambda$ we say that $f_{\mu}$ has uniformly development at $x=b(\mu)$
\[
f_{\mu}(x)\sim_{b(\mu)} \sum_{i=1}^n a_i(\mu)(b(\mu)-x)^{-\alpha_i(\mu)} \text{ in }\mu_0
\]
if there exists an open neighbourhood $U$ of $\mu_0$ such that for all $\hat\mu\in U$,
\[
\lim_{(x,\mu)\rightarrow(b(\hat\mu),\hat\mu)} (b(\mu)-x)^{\alpha_k(\mu)}\left(f_{\mu}(x)-\sum_{i=1}^{k-1} a_i(\mu)(b(\mu)-x)^{-\alpha_i(\mu)}\right)=a_k(\hat\mu)\neq 0
\]
for each $k=1,\dots,n$. If $b\equiv +\infty$ then 
\[
f_{\mu}(x)\sim_{+\infty} \sum_{i=1}^n a_i(\mu)x^{\alpha_i(\mu)} \text{ in }\mu_0
\]
if for all $\hat\mu\in U$,
\[
\lim_{(x,\mu)\rightarrow(b(\hat\mu),\hat\mu)} x^{-\alpha_k(\mu)}\left(f_{\mu}(x)-\sum_{i=1}^{k-1} a_i(\mu)x^{\alpha_i(\mu)}\right)=a_k(\hat\mu)\neq 0.
\]
We use the analogous definition at $a$.
\end{defi}

Similarly as in Definition~\ref{def:cont_quantifiable}, the functions $\mu\mapsto \alpha_k(\mu)$ are continuous at $\mu=\mu_0$.

\begin{defi}\label{def:fm}
Let $f\in C[0,+\infty)$. Setting $[f]_0\!:=f$, we define
\[
[f]_m(x)\!:=x^2[f]_{m-1}(x)+x\int_0^x [f]_{m-1}(s)ds
\]
for all $m\geq 1$.
\end{defi}

The functions $\F[f]$ and $\F\bigl[[f]_m\bigr]$ are related by the following result.

\begin{lema}[see \cite{ManRojVil2016}]\label{lema:F}
Let $f\in C[0,+\infty)$. For any $m\geq 1$,
\[
\F[f](x)=\frac{1}{x^{2m}}\F\bigl[[f]_m\bigr](x) \text{ for all }x>0.
\]
\end{lema}

In the following statement the assumption $M[[L_{\mu}]_0]\equiv \cdots\equiv M[[L_{\mu}]_{\ell-2}]\equiv 0$ in assertion $(a)$ is void in case that $\ell=1$.

\begin{prop}\label{prop:fm}
Let $\{f_{\mu}\}_{\mu\in\Lambda}$ be a continuous family of continuous functions on $[0,+\infty)$ satisfying $f_{\mu}(x)\sim_{\infty} \sum_{i=1}^N a_i(\mu) x^{-2n_i} + b(\mu)x^{\beta(\mu)}$ at $\mu_0$ with $1\leq n_1 <n_2<\cdots<n_N$ positive integers and $\beta(\mu_0)<-2n_N$. The following holds:
\begin{enumerate}[$(a)$]
\item If $M[[f_{\mu}]_0]\equiv \cdots \equiv M[[f_{\mu}]_{\ell-2}]\equiv 0$ and $M[[f_{\mu_0}]_{\ell-1}]\neq 0$ with $1\leq\ell\leq n_N$ then, for $m=1,2,\dots,\ell-1,$
\[
[f_{\mu}]_m(x)\sim_{\infty} \sum_{\begin{smallmatrix}i=1\\ n_i>m\end{smallmatrix} }^N \prod_{j=1}^{m} \frac{2j-2n_i}{2j-2n_i-1}a_i(\mu) x^{2m-2n_i}+\prod_{j=1}^{m} \frac{\beta(\mu)+2j}{\beta(\mu)+2j-1}b(\mu)x^{\beta(\mu)+2m} \text{ at }\mu_0
\]
and $[f_{\mu}]_{\ell}(x)\sim_{\infty} M[[f_{\mu}]_{\ell-1}]x$ at $\mu_0$.
\item If $M[[f_{\mu}]_0]\equiv \cdots \equiv M[[f_{\mu}]_{n_N-1}]\equiv 0$ then 
$[f_{\mu}]_{n_N}(x)\sim_{\infty} \prod_{j=1}^{n_N} \frac{\beta(\mu)+2j}{\beta(\mu)+2j-1}b(\mu) x^{\beta(\mu)+2n_N}$ at $\mu_0$.
\end{enumerate}
\end{prop}

\begin{prova}
We shall prove the result by induction on $m$. To do so, we shall first assume $m\in\{1,\dots,n_1-1\}$. (In the case $n_1=1$ this assumption is void and we move to the next step.) Let us start considering $m=1$. By definition of $[f]_1$ and elementary manipulations, we have
\begin{align*}
\frac{[f_{\mu}]_{1}(x)-\sum_{i=1}^{k-1}\frac{2-2n_i}{1-2n_i}a_i(\mu)x^{2-2n_i}}{x^{2-2n_k}}&=\frac{x^2[f_{\mu}]_0(x)-\sum_{i=1}^{k-1}\frac{2-2n_i}{1-2n_i}a_i(\mu)x^{2-2n_i}+x\int_0^x [f_{\mu}]_0(s)ds}{x^{2-2n_k}}\\
&= \frac{[f_{\mu}]_0(x)-\sum_{i=1}^{k-1}a_i(\mu)x^{-2n_i}}{x^{-2n_k}}+ \frac{\int_0^x [f_{\mu}]_0(s)ds-\sum_{i=1}^{k-1}\frac{a_i(\mu)}{1-2n_i}x^{1-2n_i} }{x^{1-2n_k}}
\end{align*}
for any $k=1,\dots,N$. The hypothesis $M[[f_{\mu}]_0]\equiv 0$ together with $n_i\geq 1$ imply that the numerator of the second quotient tends to zero as $x$ tends to infinity. Therefore the second quotient is a $0/0$-indeterminacy as $(x,\mu)\rightarrow(+\infty,\hat\mu)$ for any $\hat\mu\approx\mu_0$. We apply the Uniform H\^opital's Rule in \cite[Proposition A.1]{ManRojVil2016} to deduce that
\[
\lim_{(x,\mu)\rightarrow(+\infty,\hat\mu)} \frac{[f_{\mu}]_1(x)-\sum_{i=1}^{k-1}\frac{2-2n_i}{1-2n_i}a_i(\mu)x^{2-2n_i}}{x^{2-2n_k}}=\frac{2-2n_k}{1-2n_k}a_k(\hat\mu)
\]
for any $k=1,\dots,N$. Similarly,
\[
\lim_{(x,\mu)\rightarrow(+\infty,\hat\mu)} \frac{[f_{\mu}]_1(x)-\sum_{i=1}^{N}\frac{2-2n_i}{1-2n_i}a_i(\mu)x^{2-2n_i}}{x^{\beta(\mu)+2}}=\frac{\beta(\mu)+2}{\beta(\mu)+1}b(\hat\mu).
\]
Therefore,
\[
[f_{\mu}]_1(x)\sim_{\infty}  \sum_{i=1}^{N}\frac{2-2n_i}{1-2n_i}a_i(\mu)x^{2-2n_i} + \frac{\beta(\mu)+2}{\beta(\mu)+1}b(\hat\mu) x^{\beta(\mu)+2} \text{ at }\mu_0.
\]
This proves the result for $m=1$. If $m\in\{2,\dots,n_1-1\}$ the result follows identically changing $[f_{\mu}]_0$ by $[f_{\mu}]_1$. So we have
\[
[f_{\mu}]_m(x)\sim_{\infty} \sum_{i=1}^{N}\prod_{j=1}^m \frac{2j-2n_i}{2j-2n_i-1}a_i(\mu)x^{2m-2n_i} + \prod_{j=1}^m \frac{\beta(\mu)+2j}{\beta(\mu)+2j-1}b(\mu)x^{\beta(\mu)+2m} \text{ at }\mu_0
\]
for $m=1,\dots,n_1-1$. Let us consider now that $m\in\{n_1,\dots,n_2-1\}$, and let us start by taking $m=n_1$. The same procedure as before can be applied taking into account that the fist term in the asymptotic development disappear. Indeed, we have
\[
[f_{\mu}]_{n_1}(x)\sim_{\infty} \sum_{i=2}^{N}\prod_{j=1}^{n_1} \frac{2j-2n_i}{2j-2n_i-1}a_i(\mu)x^{2n_1-2n_i} + \prod_{j=1}^{n_1} \frac{\beta(\mu)+2j}{\beta(\mu)+2j-1}b(\mu)x^{\beta(\mu)+2n_1} \text{ at }\mu_0,
\]
where now the sum starts at $i=2$ instead of the previous $i=1$. Now using the first procedure the result follows for $m=n_1,\dots,n_2-1$. Inductively the result holds for $m=1,\dots,\ell-1$ provided that $M[[f_{\mu}]_0]\equiv \cdots M[[f_{\mu}]_{\ell-2}]\equiv 0$. To finish the proof of $(a)$, let us assume that $M[[f_{\mu_0}]_{\ell-1}]\neq 0$. Then
\[
\lim_{(x,\mu)\rightarrow(+\infty,\hat\mu)}\frac{[f_{\mu}]_{\ell}(x)}{x}=\lim_{(x,\mu)\rightarrow(+\infty,\hat\mu)}\left(x[f_{\mu}]_{\ell-1}(x)+\int_0^x [f_{\mu}]_{\ell-1}(s)ds\right) = M[[f_{\hat\mu}]_{\ell-1}]\neq 0
\]
for all $\hat\mu\approx\mu_0$. Here we used that 
\[
[f_{\mu}]_{\ell-1}(x)\sim_{\infty} \sum_{\begin{smallmatrix}i=1\\ n_i>\ell-1\end{smallmatrix} }^N \prod_{j=1}^{\ell-1} \frac{2j-2n_i}{2j-2n_i-1}a_i(\mu)x^{2\ell-2-2n_i}
\]
and, since $\ell\leq n_i$, $x [f_{\mu}]_{\ell-1}(x) \rightarrow 0$ as $(x,\mu)\rightarrow(+\infty,\hat\mu)$. This proves $(a)$. To show $(b)$ let us assume that $M[[f_{\mu}]_0]\equiv \cdots M[[f_{\mu}]_{n_N-1}]\equiv 0$. The procedure holds for all $m=1,\dots,n_N$ in this case. For $m=n_N$ all the terms in the asymptotic development associated with even powers disappear and so the only remaining term is the one with betas. Then,
\[
[f_{\mu}]_{n_N}(x)\sim_{\infty} \prod_{j=1}^{n_N} \frac{\beta(\mu)+2j}{\beta(\mu)+2j-1}b(\mu)x^{\beta(\mu)+2n_N} \text{ at }\mu_0.
\] 
This ends the proof of the result.
\end{prova}

\begin{rem}
A very useful tool for the computation of momenta was introduced in~\cite{ManRojVil2016}. If $f$ is quantifiable by $\alpha<-2n+1$, $n\in\N$ with $n\geq 2$, and $M[[f]_i]=0$ for $i=0,\dots,n-2$ then
\[
M[[f]_{n-1}]= C_n \int_0^{\infty} x^{2n-2}f(x)dx = C_n M_{n}[f]
\]
for some constant $C_n\neq 0$. We point out that this result do not contemplate the case when $\alpha=-2n$ and $M[[f]_i]=0$ for $i=0,\dots,n-2$. Therefore this simplification can not be used in Proposition~\ref{prop:fm}.
\end{rem}

From now on, for the sake of simplicity, we shall assume that the functions are analytic on $[0,+\infty)$. The reason is that in Section~\ref{sec:dynamic} the differential system~\eqref{centro} is assumed to be analytic and so the functions involved also are. However, the reader may notice that weaker regularity is allowed in the forthcoming definitions and results. Let us recall at this point the notions of Chebyshev system and its relation with the Wronskian, which are both key ingredients for our purposes. 

\begin{defi}\label{defi_ECT}
Let $f_0,f_1,\dots f_{n-1}$ be analytic functions on an open real interval $I$. The ordered set
$(f_0,f_1,\dots f_{n-1})$ is an \emph{extended complete Chebyshev system} (for short, a ECT-system) on $I$ if, for all $k=1,2,\dots n$, any nontrivial linear combination
\[
\alpha_0f_0(x)+\alpha_1f_1(x)+\cdots+\alpha_{k-1}f_{k-1}(x)
\]
has at most $k-1$ isolated zeros on $I$ counted with multiplicities.
(Let us mention that, in these abbreviations, ``T'' stands for Tchebycheff, which in some sources is the transcription of the Russian name Chebyshev).
\end{defi}

\begin{defi}\label{wronskian}
Let $f_0,f_1,\dots,f_{n-1}$ be analytic functions on an open interval $I$ of $\R$. Then
\[
W[f_0,f_1,\dots,f_{n-1}](x)=\text{det}\left(f_j^{(i)}(x)\right)_{0\leq i,j\leq n-1}=
\left|\begin{array}{ccc}
f_0(x) & \cdots & f_{n-1}(x)\\
f_0'(x) & \cdots & f_{n-1}'(x)\\
		&	\vdots & \\
f_0^{(n-1)}(x) & \cdots & f_{n-1}^{(n-1)}(x)
\end{array}\right|
\]
is the \emph{Wronskian} of $(f_0,f_1,\dots,f_{n-1})$ at $x\in I$.
\end{defi}

These two notions are closely related by the following result (see for instance \cite{Karlin}).

\begin{lema}\label{lema:ECT-Wronskia}
$(f_0,f_1,\dots,f_{n-1})$ is an ECT-system on $I$ if and only if, for each $k=1,2,\dots,n$,
\[
W[f_0,f_1,\dots,f_{k-1}](x)\neq 0\ \ \text{for all }x\in I.
\]
\end{lema}

In the study of bifurcation of critical periodic orbits from the outer boundary the main objective is to bound the zeros of the derivative of the period function for energy levels $h\approx h_0(\mu)$ uniformly on the parameters close to a fixed $\mu_0\in\Lambda$. Invoking equality~\eqref{FT}, this problem is tackled by completing $\F[f_{\mu}]$ with a collection of analytic functions $\phi_{\mu}^1,\dots,\phi_{\mu}^n$ in order that $(\phi_{\mu}^1,\dots,\phi_{\mu}^n,\F[f_{\mu}])$ form an ECT-system on $(h_0(\mu)-\epsilon,h_0(\mu))$ for some $\epsilon>0$ and all $\mu\approx\mu_0$. Note that guarantee the uniformity with respect to the parameters of the system is mandatory to obtain the desired upper bounds. With this aim in view, and on account of the characterization in Lemma~\ref{lema:ECT-Wronskia}, given $\nu_1,\nu_2,\dots,\nu_n\in\R$, we consider the linear ordinary differential operator
\[
\LL_{\boldsymbol\nu_n(\mu)}:C^{\omega}(0,+\infty)\rightarrow C^{\omega}(0,+\infty)
\]
defined by
\[
\LL_{\boldsymbol\nu_n(\mu)}[f](x)\!:= \frac{W[x^{\nu_1},x^{\nu_2},\dots,x^{\nu_n},f(x)]}{x^{\sum_{i=1}^n(\nu_i-i)}}.
\]
Here, and in what follows, for the sake of shortness we use the notation $\boldsymbol\nu_n=(\nu_1,\nu_2,\dots,\nu_n)$. Furthermore we define $\LL_{\boldsymbol\nu_0}=id$ in order that forthcoming statements contemplate the case $n=0$ as well. We also denote by $C^{\omega}(0,+\infty)$ the analytic functions on $(0,+\infty)$ and by $C^{\omega}[0,+\infty)$ the functions on $C^{\omega}(0,+\infty)$ that can be extended analytically to $x=0$.

\begin{prop}[See~\cite{ManRojVil2016a}]\label{prop:commutar}
For any $f\in C^{\omega}(0,+\infty)$ and $\nu_1,\dots,\nu_n\in\R$, the following recurrence holds:
\[
\LL_{\boldsymbol\nu_n}[f](x)=c_n\left( x \LL_{\boldsymbol\nu_{n-1}}[f]'(x)-\nu_n \LL_{\boldsymbol\nu_{n-1}}[f](x)\right),
\]
where $c_1\!:=1$ and $c_n\!:=\prod_{i=1}^{n-1}(\nu_n-\nu_i)$ for $n\geq 2$. In particular, if $f\in C^{\omega}[0,+\infty)$, then $\LL_{\boldsymbol\nu_n}[f]\in C^{\omega}[0,+\infty)$. Moreover, $\F\circ \LL_{\boldsymbol\nu_n}=\LL_{\boldsymbol\nu_n}\circ \F$.
\end{prop}

In the statement of the following result the assumptions $M[[L_{\mu}]_0]\equiv \cdots\equiv M[[L_{\mu}]_{\ell-2}]\equiv 0$ in assertion $(a)$ and $M_i[[L_{\mu}]_{n_N}]\equiv 0$ for $i=1,\dots,\ell-1$ in assertion $(d1)$ are void in case that $\ell=1$.

\begin{thmx}\label{thm:tecnic}
Let $\Lambda$ be an open subset of $\R^d$ and $\{f_{\mu}\}_{\mu\in\Lambda}$ be a continuous family of analytic functions on $[0,+\infty)$. Assume that, in a neighbourhood of some fixed $\mu_0\in\Lambda$ there exist $n\geq 0$ continuous functions $\nu_1,\nu_2,\dots,\nu_n$ pairwise distinct at $\mu=\mu_0$ such that the function $L_{\mu}:=\LL_{\boldsymbol\nu_n(\mu)}[f_{\mu}]$ satisfies
\[
L_{\mu}(x)\sim_{\infty} \sum_{i=1}^N a_i(\mu) x^{-2n_i} + b(\mu)x^{\beta(\mu)}\text{ at }\mu_0
\]
with $1\leq n_1<n_2<\dots<n_N$ positive integers and $\beta(\mu_0)<-2n_N$. The following holds:
\begin{enumerate}[$(a)$]
\item If $M[[L_{\mu}]_0]\equiv \cdots\equiv M[[L_{\mu}]_{\ell-2}]\equiv 0$ and $M[[L_{\mu_0}]_{\ell-1}]\neq 0$ for some $1\leq\ell\leq n_N$ then 
\[
(\LL_{\boldsymbol\nu_n(\mu)}\circ\F)[f_{\mu}](x)\sim_{\infty} M[[L_{\mu}]_{\ell-1}]x^{1-2\ell} \text{ at }\mu_0.
\]
\item If $M[[L_{\mu}]_0]\equiv \cdots\equiv M[[L_{\mu}]_{n_N-1}]\equiv 0$ and $\beta(\mu_0)+2n_N=-1$ then
\[
(\LL_{\boldsymbol\nu_n(\mu)}\circ\F)[f_{\mu}](x)\sim_{\infty} \prod_{j=1}^{n_N} \frac{\beta(\mu)+2j}{\beta(\mu)+2j-1}b(\mu) \frac{\Omega(x,\beta(\mu)+2n_N)}{x^{1+2n_N}} \text{ at }\mu_0.
\]
\item If $M[[L_{\mu}]_0]\equiv \cdots\equiv M[[L_{\mu}]_{n_N-1}]\equiv 0$ and $\beta(\mu_0)+2n_N>-1$ then
\[
(\LL_{\boldsymbol\nu_n(\mu)}\circ\F)[f_{\mu}](x)\sim_{\infty} \fun\!\left(\beta(\mu)+2n_N\right)\prod_{j=1}^{n_N} \frac{\beta(\mu)+2j}{\beta(\mu)+2j-1}b(\mu)x^{\beta(\mu)} \text{ at }\mu_0.
\]
\item If $M[[L_{\mu}]_0]\equiv \cdots\equiv M[[L_{\mu}]_{n_N-1}]\equiv 0$ and $\beta(\mu_0)+2n_N<-1$, let us take $m\in\N$ such that $\beta(\mu_0)+2n_N+2m\in[-1,1)$ and let us assume additionally that $\beta(\mu_0)+2n_N+2m\neq 0$. In this case,
\begin{enumerate}[$(d1)$]
\item If $M_i[[L_{\mu}]_{n_N}]\equiv 0$ for $i=1,\dots,\ell-1$ and $M_{\ell}[[L_{\mu_0}]_{n_N}]\neq 0$ for some $1\leq \ell\leq m$, then 
\[
(\LL_{\boldsymbol\nu_n(\mu)}\circ\F)[f_{\mu}](x)\sim_{\infty} M_{\ell}[[L_{\mu}]_{n_N}]x^{1-2n_N-2\ell} \text{ at }\mu_0.
\]
\item If $M_i[[L_{\mu}]_{n_N}]\equiv 0$ for $i=1,\dots,m$ and $\beta(\mu_0)+2n_N+2m\neq -1$ then
\[
(\LL_{\boldsymbol\nu_n(\mu)}\circ\F)[f_{\mu}](x)\sim_{\infty} \prod_{j=1}^{n_N+m} \frac{\beta(\mu)+2j}{\beta(\mu)+2j-1}b(\mu)\fun\!(\beta(\mu)+2n_N+2m) x^{\beta(\mu)} \text{ at }\mu_0.
\]
\item If $M_i[[L_{\mu}]_{n_N}]\equiv 0$ for $i=1,\dots,m$ and $\beta(\mu_0)+2n_N+2m=-1$ then 
\[
(\LL_{\boldsymbol\nu_n(\mu)}\circ\F)[f_{\mu}](x)\sim_{\infty} \prod_{j=1}^{n_N+m} \frac{\beta(\mu)+2j}{\beta(\mu)+2j-1}b(\mu)\frac{\Omega(x,\beta(\mu)+2n_N+2m)}{x^{2n_N+2m+1}} \text{ at }\mu_0.
\]
\end{enumerate}
\end{enumerate}
\end{thmx}

\begin{prova}
Let us start by considering $n=0$. In this case $L_{\mu}=f_{\mu}$ and the assumptions of Proposition~\ref{prop:fm} are satisfied. Moreover, if $(a)$ is satisfied then by $(a)$ in Proposition~\ref{prop:fm} we have
\[
[f_{\mu}]_{\ell}(x)\sim_{\infty} M[[f_{\mu}]_{\ell-1}]x \text{ at }\mu_0.
\]
Therefore, using Lemma~\ref{lema:F} and Theorem~\ref{thm:quantificar} with $[f_{\mu}]_{\ell}$, since $\fun(1)=1$,
\[
\F[f_{\mu}](x)=\frac{1}{x^{2\ell}}\F[[f_{\mu}]_{\ell}](x)\sim_{+\infty} M[[f_{\mu}]_{\ell-1}]x^{1-2\ell} \text{ at }\mu_0.
\]
This proves $(a)$. Let us show $(b)$, $(c)$ and $(d)$. By hypothesis we have $M[[f_{\mu}]_i]\equiv 0$ for $i=0,\dots,n_N-1$. Therefore by $(b)$ in Proposition~\ref{prop:fm},
\[
[f_{\mu}]_{n_N}\sim_{+\infty}\prod_{j=1}^{n_N} \frac{\beta(\mu)+2j}{\beta(\mu)+2j-1}b(\mu) x^{\beta(\mu)+2n_N} \text{ at }\mu_0.
\]
The result follows using Theorem~\ref{thm:quantificar} with $[f_{\mu}]_{n_N}$ and Lemma~\ref{lema:F} again. This ends the proof for the case $n=0$.

Let us consider now $n\geq 1$. By Proposition~\ref{prop:commutar} $L_{\mu}=\LL_{\boldsymbol\nu_n(\mu)}[f_{\mu}]$ is an analytic function on $[0,+\infty)$ for each $\mu\in\Lambda$ and
\[
(\LL_{\boldsymbol\nu_n(\mu)}\circ \F)[f_{\mu}](x)=(\F\circ \LL_{\boldsymbol\nu_n(\mu)})[f_{\mu}](x).
\]
Then the result follows by applying the case $n=0$ to the family $\{\LL_{\boldsymbol\nu_n(\mu)}[f_{\mu}]\}_{\mu\in\Lambda}$.
\end{prova}

\begin{rem}
Let us recover at this point the examples in~\eqref{examples}. The function $f$ satisfies the assumptions of Theorem~\ref{thm:tecnic} with $n=0$, $n_1=1$ and $\beta=-5/2$. Since $\beta+2=\frac{-1}{2}>-1$ then assertion $(b)$ of the Theorem states that
\[
\F[f](x)\sim_{+\infty} \frac{\sqrt{\pi}~\Gamma\!\left(\frac{1}{4}\right)}{6~\Gamma\!\left(\frac{3}{4}\right)}x^{-\frac{5}{2}}.
\]
The function $g$ satisfies the hypothesis with $n=0$, $n_1=1$ and $\beta=-5$. In this case $\beta+2=-3<-1$ and
\[
M[[g]_1]=\int_0^{\infty}\left(\begin{cases}
\frac{3}{4x^3} & \text{ if }x\geq 1\\
\frac{3}{4}x^2(13x-12) & \text{ if }x\in[0,1)
\end{cases}\right)dx=-\frac{3}{16}.
\]
That is, $g$ satisfies the hypothesis in assertion $(c1)$ of the Theorem with $j=1$. Therefore,
\[
\F[g](x)\sim_{+\infty} -\frac{3}{16}x^{-3}.
\]
Finally, the function $h$ satisfies the hypothesis with $n=0$, $n_1=1$ and $\beta=-3$. In this example assumptions in assertion~$(b)$ are satisfied and
\[
\F[h](x)\sim_{+\infty} \frac{\log x}{2x^3}.
\]
\end{rem}

We point out that Theorem~\ref{thm:tecnic} together with Theorem~\ref{thm:quantificar} cover all possible situations of the uniform asymptotic development of $f_{\mu}$ except the case when all powers are negative even numbers at $\mu=\mu_0$. The last result of this Section is a recursive formula for the computation of a certain momentum that will be useful in the applications.

\begin{lema}\label{lema:calcul-moment}
Let $f\in C^{\omega}[0,+\infty)$, $n\geq 1$ and $\nu_1,\nu_2,\dots,\nu_n\in\R$. Then,
\[
M\bigl[[\LL_{\boldsymbol\nu_n}[f]]_1\bigr]=c_n\lim_{R\rightarrow+\infty}\!\left( R^3\LL_{\boldsymbol\nu_{n-1}}[f](R) + \frac{(\nu_n+1)R^2}{2}\!\!\int_0^R \!\!\!\LL_{\boldsymbol\nu_{n-1}}[f](x)dx - \frac{3\nu_n+5}{2}\!\!\int_0^R \!\!\!x^2\LL_{\boldsymbol\nu_{n-1}}[f](x)dx  \right)
\]
where $c_1\!:=1$ and $c_n\!:=\prod_{i=1}^{n-1}(\nu_n-\nu_i)$ for $n\geq 2$.
\end{lema}

\begin{prova}
By Definition~\ref{def:fm} we have
\[
M\bigl[[\LL_{\boldsymbol\nu_n}[f]]_1\bigr]=\int_0^{\infty}[\LL_{\boldsymbol\nu_n}[f]]_1(x)dx=\lim_{R\rightarrow+\infty}\int_{0}^R \left(x^2\LL_{\boldsymbol\nu_n}[f](x) + x\int_0^x \LL_{\boldsymbol\nu_n}[f](s)ds\right) dx.
\]
Using the recursive expression in Proposition~\ref{prop:commutar},
\begin{align*}
M\bigl[[\LL_{\boldsymbol\nu_n}[f]]_1\bigr]=c_n\lim_{R\rightarrow+\infty}&\left(\int_{0}^R \bigl(x^3\LL_{\boldsymbol\nu_{n-1}}[f]'(x)-\nu_n x^2\LL_{\boldsymbol\nu_{n-1}}[f](x)\bigr)dx\right.+\\
&+\left.\int_0^R \left(x\int_0^x \bigl(s\LL_{\boldsymbol\nu_{n-1}}[f]'(s)-\nu_n\LL_{\boldsymbol\nu_{n-1}}[f](s)\bigr)ds\right) dx\right).
\end{align*}
The result then follows integrating by parts.
\end{prova}

\section{Criticality of the period function at the outer boundary}\label{sec:dynamic}

In this section we apply Theorem~\ref{thm:tecnic} in order to obtain sufficient conditions to bound the number of critical periodic orbits that may bifurcate from the outer boundary of the period annulus in families of planar potential centers. We consider analytic differential systems~\eqref{centro} depending on a parameter $\mu\in\Lambda\subset\R^d$ and we assume that the origin is a non-degenerate center for all $\mu$. We denote by $\mathcal I_{\mu}=(x_{\ell}(\mu),x_r(\mu))$ the projection of the period annulus on the $x$-axis, $x_{\ell}<0<x_r$ and by $h_0(\mu)$ the energy level at the outer boundary of the period annulus.

\begin{defi}\label{def:Ak}
We say that the family~\eqref{centro} verifies the hypothesis \textnormal{\textbf{(H)}} in case that:
\begin{enumerate}[$(a)$]
\item For all $k\geq 0$, the map $(x,\mu)\longmapsto V_{\mu}^{(k)}(x)$ is continuous on $\{(x,\mu)\in\R\times\Lambda : x\in I_{\mu}\},$
\item $\mu\longmapsto x_{r}(\mu)$ is continuous on $\Lambda$ or $x_{r}(\mu)=+\infty$ for all $\mu\in\Lambda,$
\item $\mu\longmapsto x_{\ell}(\mu)$ is continuous on $\Lambda$ or $x_{\ell}(\mu)=-\infty$ for all $\mu\in\Lambda,$
\item $\mu\longmapsto h_0(\mu)$ is continuous on $\Lambda$ or $h_0(\mu)=+\infty$ for all $\mu\in\Lambda.$
\end{enumerate}\vspace*{-.5cm}
\end{defi}

\begin{lema}[See~\cite{ManRojVil2016}]\label{lema:g_menys_1_continua}
Let $\{X_{\mu}\}_{\mu\in\Lambda}$ be a  family of potential analytic differential systems verifying \textnormal{\textbf{(H)}}.  
Then the map $(z,\mu)\longmapsto g_{\mu}^{-1}(z)$ is continuous on the open set
$\bigl\{(z,\mu)\in\R\times\Lambda: z\in
\bigl(-\sqrt{h_{0}(\mu)},\sqrt{h_{0}(\mu)}\bigr)\bigr\}$.
\end{lema}

This section is divided in two parts according with the dichotomy produced by $h_0$. First, Section~\ref{sec:infinite} deals with the case $h_0\equiv+\infty$. Second, Section~\ref{sec:finite} is dedicated to the case $h_0$ finite.

\subsection{Potential systems with infinite energy}\label{sec:infinite}

In this section we present sufficient conditions to bound the criticality at the outer boundary for potential systems satisfying $h_0(\mu)=+\infty$ for all $\mu\in\Lambda$. Following the strategy in~\cite{ManRojVil2016,ManRojVil2016a,Rojas2018}, we find sufficient conditions such that $fT_{\mu}'$ can be embedded into the ECT-system $(h^{\nu_1(\mu)},h^{\nu_2(\mu)},\dots,h^{\nu_n(\mu)})$, where $f$ is an analytic non-vanishing function. Next result is a combination of~\cite[Lemma~3.5]{ManRojVil2016a} and~\cite[Lemma~3.3]{Rojas2018}, and it is the key piece that connects the analytic tools studied in Section~\ref{sec:technical} with the dynamical results we are looking for.

\begin{lema}\label{lema:crit}
Let $\{X_{\mu}\}_{\mu\in\Lambda}$ be a family of potential analytic differential systems verifying \textnormal{\textbf{(H)}} and such that $h_0\equiv +\infty$. Assume that there exist $n\geq 1$ continuous functions $\nu_1,\nu_2,\dots,\nu_n$ in a neighbourhood of some fixed $\mu_0\in\Lambda$, a continuous function $\alpha:\Lambda\rightarrow\R$ with $\alpha(\mu_0)=-1$ and an analytic non-vanishing function $f$ on $(0,+\infty)$ such that
\[
\lim_{(h,\mu)\rightarrow(+\infty,\mu_0)}\frac{h^{\nu_n(\mu)}}{\Omega(h,\alpha(\mu))^m}W[h^{\nu_1(\mu)},\dots,h^{\nu_{n-1}(\mu)},f(h)T_{\mu}'(h)]=\ell\neq 0
\]
with $m\in\{0,1\}$. Then $\mathrm{Crit}\bigl((\out_{\mu_0},X_{\mu_0}),X_{\mu}\bigr)\leqslant n-1$.
\end{lema}

The assumption requiring the existence of functions $\nu_1,\nu_2,\dots,\nu_n$ in the following statement is void in case that $n=0$. The same happens to the assumptions $M[[L_{\mu}]_0]\equiv \cdots\equiv M[[L_{\mu}]_{\ell-2}]\equiv 0$ in assertion $(a)$ and $M_i[[L_{\mu}]_{n_N}]\equiv 0$ for $i=1,\dots,\ell-1$ in assertion $(c1)$ in case that $\ell=1$.

\begin{thmx}\label{thm:criticalitat-infinite}
Let $\{X_{\mu}\}_{\mu\in\Lambda}$ be a family of potential analytic differential systems verifying \textnormal{\textbf{(H)}} with $h_0\equiv+\infty$ and that there exist $n\geq 0$ continuous functions $\nu_1,\nu_2,\ldots,\nu_n$ in a neighbourhood of some fixed $\mu_0\in\Lambda$ such that 
\[
L_{\mu}(x)\!:=\LL_{\boldsymbol\nu_n(\mu)}[x(g_{\mu}^{-1})''(x)-x(g_{\mu}^{-1})''(-x)]\sim_{+\infty} \sum_{i=1}^N a_i(\mu) x^{-2n_i} + b(\mu)x^{\beta(\mu)}\text{ at }\mu_0
\]
with $1\leq n_1<n_2<\dots<n_N$ positive integers and $\beta(\mu_0)<-2n_N$. 
Then $\mathrm{Crit}\bigl((\out_{\mu_0},X_{\mu_0}),X_{\mu}\bigr)\leqslant n$ if one of the following assertions hold:
\begin{enumerate}[$(a)$]
\item If $M[[L_{\mu}]_0]\equiv \cdots\equiv M[[L_{\mu}]_{\ell-2}]\equiv 0$ and $M[[L_{\mu_0}]_{\ell-1}]\neq 0$ for some $1\leq\ell\leq n_N$.
\item If $M[[L_{\mu}]_0]\equiv \cdots\equiv M[[L_{\mu}]_{n_N-1}]\equiv 0$ and $\beta(\mu_0)+2n_N\geq -1$.
\item If $M[[L_{\mu}]_0]\equiv \cdots\equiv M[[L_{\mu}]_{n_N-1}]\equiv 0$ and $\beta(\mu_0)+2n_N<-1$, let us take $m\in\N$ such that $\beta(\mu_0)+2n_N+2m\in[-1,1)\setminus\{0\}$. In this case:
\begin{enumerate}[$(c1)$]
\item If $M_i[[L_{\mu}]_{n_N}]\equiv 0$ for $i=1,\dots,\ell-1$ and $M_{\ell}[[L_{\mu_0}]_{n_N}]\neq 0$ for some $1\leq \ell\leq m$.
\item If $M_i[[L_{\mu}]_{n_N}]\equiv 0$ for $i=1,\dots,m$.
\end{enumerate}
\end{enumerate}
\end{thmx}

\begin{prova}
For the sake of shortness let us denote $f_{\mu}(x)\!:=x(g_{\mu}^{-1})''(x)-x(g_{\mu}^{-1})''(-x)$. Lemma~\ref{lema:g_menys_1_continua} and the hypothesis~\textnormal{\textbf{(H)}} imply that $\{f_{\mu}\}_{\mu\in\Lambda}$ is a continuous family of analytic functions on $[0,+\infty)$. According to equality~\eqref{FT} the result will follow once we show that there exist $M,\epsilon>0$ in such a way $\F[f_{\mu}]$ has at most $n$ isolated zeros for $h>M$ and $\|\mu-\mu_0\|<\epsilon$, multiplicities taken into account.

Let us assume that $L_{\mu}\!:=\LL_{\boldsymbol\nu_n(\mu)}[f_{\mu}]$ satisfies one of the hypothesis of the statement. Therefore $L_{\mu}$ satisfies one of the hypothesis in Theorem~\ref{thm:tecnic}, so we can assert that either
\[
(\LL_{\boldsymbol\nu_n(\mu)}\circ\F)[f_{\mu}](h)\sim_{+\infty} C(\mu) h^{\xi(\mu)} \text{ at } \mu_0
\]
or 
\[
(\LL_{\boldsymbol\nu_n(\mu)}\circ\F)[f_{\mu}](h)\sim_{+\infty} C(\mu) \Omega(h,\alpha(\mu))h^{\xi(\mu)} \text{ at } \mu_0
\]
for some functions $C,\xi,\alpha$ with $C(\mu_0)\neq 0$ and $\alpha(\mu_0)=-1$. Taking into account the definition of the operator $\LL_{\boldsymbol\nu_n(\mu)}$, we have that either
\[
\lim_{(x,\mu)\rightarrow(+\infty,\mu_0)}\frac{W\bigr[h^{\nu_1(\mu)},\dots,h^{\nu_n(\mu)},\F[f_{\mu}](h)\bigr]}{h^{\xi(\mu)+\sum_{i=1}^{n}(\nu_i(\mu)-i)}}= C(\mu_0),
\]
or 
\[
\lim_{(x,\mu)\rightarrow(+\infty,\mu_0)}\frac{W\bigr[h^{\nu_1(\mu)},\dots,h^{\nu_n(\mu)},\F[f_{\mu}](h)\bigr]}{h^{\xi(\mu)+\sum_{i=1}^{n}(\nu_i(\mu)-i)}\Omega(h,\alpha(\mu))}= C(\mu_0).
\]
Therefore, on account of equality~\eqref{FT}, by Lemma~\ref{lema:crit} we have that $\mathrm{Crit}\bigl((\out_{\mu_0},X_{\mu_0}),X_{\mu}\bigr)\leqslant n$ as desired.
\end{prova}

\subsection{Potential systems with finite energy}\label{sec:finite}

We assume in this section that the energy at the outer boundary of system~\eqref{centro} is finite for all parameters $\mu\in\Lambda$. As in the previous works~\cite{ManRojVil2016,ManRojVil2016a,Rojas2018}, in order to embed the function $fT_{\mu}'$ into some ECT-system for an appropriate non-vanishing function $f$, the spirit of this section is to ``translate'' the case $h_0<+\infty$ to the case $h_0=+\infty$ so we can take advantage of Theorem~\ref{thm:tecnic}. This translation is given by the operator
\[
\CC\!: C^{\omega}[0,1)\longrightarrow C^{\omega}[0,+\infty)
\] 
defined by
\begin{equation}\label{defi:operator_C}
\CC[f](x)\!:=\bigl(1-\phi^2(x)\bigr)\bigl(f\circ \phi\bigr)(x)=\frac{1}{1+x^2}\bigl(f\circ \phi\bigr)(x),
\end{equation}
where $\phi(x)\!:=\frac{x}{\sqrt{1+x^2}}$. Given $\nu_1,\dots,\nu_n\in\R$, $\CC$ conjugates the operator $\LL_{\boldsymbol\nu_n}$ with the linear ordinary differential operator 
\[
\DD_{\boldsymbol\nu_n}\!: C^{\omega}(0,1)\longrightarrow C^{\omega}(0,1)
\]
defined by
\begin{equation}\label{eq:D}
\DD_{\boldsymbol\nu_n}[f](z)\!:=(z(1-z^2))^{\frac{n(n+1)}{2}}\frac{W\left[\psi_{\nu_1},\ldots,\psi_{\nu_n},f\right](z)}{\prod_{i=1}^n\psi_{\nu_i}(z)}.
\end{equation}

\begin{defi}
Let $f\in C^{\omega}[0,1).$ We call
\[
N_n[f]:=\int_0^{1} \frac{f(z)}{\sqrt{1-z^2}}\left(\frac{z}{\sqrt{1-z^2}}\right)^{2n-2}dz
\]
the $n$-th \emph{momentum} of $f$, whenever it is well defined. If $n=1$ we simply say that $N[f]\!:=N_1[f]$ is the momentum of $f$.
\end{defi}

\begin{lema}\label{lem:lemes_varis}
Consider $\nu_1,\nu_2,\dots,\nu_n\in\R$. Then the following hold:
\begin{enumerate}[$(a)$]
\item $\CC[\psi_{\nu_i}](x)=x^{\nu_i}$ for $i=1,2,\ldots,n.$
\item $\CC\circ\DD_{\boldsymbol\nu_n}=\LL_{\boldsymbol\nu_n}\!\circ\CC$.
\item $\bigl(\F\circ \CC\bigr)[f](x)=\sqrt{1+x^2}\bigl(\CC\circ \F\bigr)[f](x)$ for any $f\in C^{\omega}(0,1).$ 
\item $N_n=M_n\circ \CC$.
\end{enumerate}
\end{lema}

\begin{defi}\label{def:moment-finit}
Given a continuous function $f$ on $[0,1)$ we define $[f]^m(z)\!:= \frac{1}{1-z^2}[\CC[f]]_{m}(\phi^{-1}(z))$.
\end{defi}

\begin{lema}\label{lema:trad}
Let $\{f_{\mu}\}_{\mu\in\Lambda}$ be a continuous family of analytic functions on $[0,1)$. Then
\[
f_{\mu}(z)\sim_{z=1} \sum_{i=1}^{n}a_i(\mu)(1-z^2)^{-\alpha_i(\mu)} \text{ at }\mu_0 \text{ if and only if } \CC[f_{\mu}](x)\sim_{+\infty} \sum_{i=1}^na_i(\mu)x^{2\alpha_i(\mu)-2} \text{ at }\mu_0.
\]
\end{lema}

\begin{prova}
By definition $\CC[f](x)=\frac{1}{1+x^2}f(\phi(x))$ with $\phi(x)=\tfrac{x}{\sqrt{1+x^2}}$. Therefore, for a fixed $\mu_0\in\Lambda$,
\[
\lim_{(x,\mu)\rightarrow(+\infty,\mu_0)}\frac{\CC[f_{\mu}](x)-\sum_{i=1}^{k-1}a_i(\mu)x^{2\alpha_i(\mu)-2}}{x^{2\alpha_k(\mu)-2}}=
\lim_{(x,\mu)\rightarrow(+\infty,\mu_0)}\frac{f_{\mu}(\phi(x))-\sum_{i=1}^{k-1}a_i(\mu)(1+x^2)^{\alpha_i(\mu)}}{(1+x^2)^{\alpha_k(\mu)}}.
\]
Denoting $z=\phi(x)$, since $\phi(x)\rightarrow 1$ as $x\rightarrow+\infty$, we have that
\[
\lim_{(x,\mu)\rightarrow(+\infty,\mu_0)}\frac{\CC[f_{\mu}](x)-\sum_{i=1}^{k-1}a_i(\mu)x^{2\alpha_i(\mu)-2}}{x^{2\alpha_k(\mu)-2}}=
\lim_{(z,\mu)\rightarrow(1,\mu_0)}\frac{f_{\mu}(z)-\sum_{i=1}^{k-1}a_i(\mu)(1-z^2)^{-\alpha_i(\mu)}}{(1-z^2)^{-\alpha_k(\mu)}}.
\]
Consequently the result follows on account of Definition~\ref{defi:asympt}.
\end{prova}

The following is an analogous version of Lemma~\ref{lema:crit} for the case $h_0$ finite. It gathers \cite[Lemma 3.10]{ManRojVil2016a} and \cite[Lemma 3.8]{Rojas2018}.

\begin{lema}\label{lema:crit-finit}
Let $\{X_{\mu}\}_{\mu\in\Lambda}$ be a family of potential analytic differential systems verifying \textnormal{\textbf{(H)}} and such that $\mu\mapsto h_0(\mu)$ is continuous on $\Lambda$. Assume that there exist $n\geq 1$ continuous functions $\nu_1,\nu_2,\dots,\nu_n$ in a neighbourhood of some fixed $\mu_0\in\Lambda$, a continuous function $\alpha:\Lambda\rightarrow\R$ with $\alpha(\mu_0)=-1$ and an analytic non-vanishing function $f$ on $(0,1)$ such that
\[
\lim_{(z,\mu)\rightarrow(1,\mu_0)}\frac{(1-z)^{\nu_n(\mu)}}{\Omega\left(\tfrac{z}{\sqrt{1-z^2}},\alpha(\mu)\right)^m}W[\psi_{\nu_1(\mu)}(z),\dots,\psi_{\nu_{n-1}(\mu)}(z),f(z)T_{\mu}'(z^2 h_0(\mu))]=\ell\neq 0
\]
with $m\in\{0,1\}$. Then $\mathrm{Crit}\bigl((\out_{\mu_0},X_{\mu_0}),X_{\mu}\bigr)\leqslant n-1$.
\end{lema}

In the same way as in the main result of the previous section, we stress that the assumption requiring the existence of functions $\nu_1,\nu_2,\dots,\nu_n$ in the following statement is void in case that $n=0$. Also are void the assumptions $N[[D_{\mu}]^0]\equiv \cdots\equiv N[[D_{\mu}]^{\ell-2}]\equiv 0$ in assertion $(a)$ and $N_i[[D_{\mu}]^{n_N}]\equiv 0$ for $i=1,\dots,\ell-1$ in assertion $(c1)$ in case that $\ell=1$.

\begin{thmx}\label{thm:criticalitat-finite}
Let $\{X_{\mu}\}_{\mu\in\Lambda}$ be a family of potential analytic differential systems verifying \textnormal{\textbf{(H)}} with $h_0(\mu)<+\infty$ for all $\mu\in\Lambda$ and that there exist $n\geq 0$ continuous functions $\nu_1,\nu_2,\ldots,\nu_n$ in a neighbourhood of some fixed $\mu_0\in\Lambda$ such that the function
\[
f_{\mu}(z)\!:=z\sqrt{h_0(\mu)}(g_{\mu}^{-1})''(z\sqrt{h_0(\mu)})-z\sqrt{h_0(\mu)}(g_{\mu}^{-1})''(-z\sqrt{h_0(\mu)})
\]
satisfies
\[
D_{\mu}(z)\!:=\DD_{\boldsymbol\nu_n(\mu)}[f_{\mu}](z)\sim_{z=1} \sum_{i=1}^N a_i(\mu) (1-z^2)^{n_i} + b(\mu)(1-z^2)^{\beta(\mu)}\text{ at }\mu_0
\]
with $0\leq n_1<n_2<\dots<n_N$ integers and $\beta(\mu_0)>n_N$.
Then $\mathrm{Crit}\bigl((\out_{\mu_0},X_{\mu_0}),X_{\mu}\bigr)\leqslant n$ if one of the following assertions hold:
\begin{enumerate}[$(a)$]
\item If $N[[D_{\mu}]^0]\equiv \cdots\equiv N[[D_{\mu}]^{\ell-2}]\equiv 0$ and $N[[D_{\mu_0}]^{\ell-1}]\neq 0$ for some $1\leq\ell\leq n_N$.
\item If $N[[D_{\mu}]^0]\equiv \cdots\equiv N[[D_{\mu}]^{n_N-1}]\equiv 0$ and $\beta(\mu_0)-n_N\leq \frac{1}{2}$.
\item If $N[[D_{\mu}]^0]\equiv \cdots\equiv N[[D_{\mu}]^{n_N-1}]\equiv 0$ and $\beta(\mu_0)-n_N>\frac{1}{2}$, let us take $m\in\N$ such that $\beta(\mu_0)-n_N-m\in(-\tfrac{1}{2},\tfrac{1}{2}]\setminus\{0\}$. In this case:
\begin{enumerate}[$(c1)$]
\item If $N_i[[D_{\mu}]^{n_N}]\equiv 0$ for $i=1,\dots,\ell-1$ and $N_{\ell}[[D_{\mu_0}]^{n_N}]\neq 0$ for some $1\leq \ell\leq m$.
\item If $N_i[[D_{\mu}]^{n_N}]\equiv 0$ for $i=1,\dots,m$.
\end{enumerate}
\end{enumerate}
\end{thmx}

\begin{prova}
According with Lemma~\ref{lema:g_menys_1_continua} and the hypothesis~\textnormal{\textbf{(H)}}, the family $\{f_{\mu}\}_{\mu\in\Lambda}$ is a continuous family of analytic functions on $[0,1)$. On account of equality~\eqref{FT}, after the appropriate rescaling,
\begin{equation}\label{FT-finit}
\F[f_{\mu}](z)=\sqrt{2}h_0(\mu) z^2T_{\mu}'(h_0(\mu)z^2) \text{ for all }z\in(0,1).
\end{equation}
Therefore the result will follow if there exist $\epsilon>0$ and a neighbourhood $U$ of $\mu_0$ such that $\F[f_{\mu}](z)$ has at most $n$ zeros for all $z\in(1-\epsilon,1)$ and $\mu\in U$, multiplicities taken into account. To this end, we shall use that the operator $\CC$ commutes the linear differential operators $\DD_{\boldsymbol\nu_n(\mu)}$, defined for functions in $C^{\omega}[0,1)$, with $\LL_{\boldsymbol\nu_n(\mu)}$, defined for functions in $C^{\omega}[0,+\infty)$. Then, as we proceeded in Theorem~\ref{thm:criticalitat-infinite}, we aim to apply Theorem~\ref{thm:tecnic} in this case to the family $\CC[f_{\mu}]$.

First, by hypothesis we have that 
\[
\DD_{\boldsymbol\nu_n(\mu)}[f_{\mu}](z)\sim_{z=1} \sum_{i=1}^N a_i(\mu) (1-z^2)^{n_i} + b(\mu)(1-z^2)^{\beta(\mu)}\text{ at }\mu_0,
\]
so, on account of Lemmas~\ref{lem:lemes_varis} and~\ref{lema:trad}, we have
\[
(\LL_{\boldsymbol\nu_n(\mu)}\circ\CC)[f_{\mu}](x)=(\CC\circ \DD_{\boldsymbol\nu_n(\mu)})[f_{\mu}](x)\sim_{+\infty} \sum_{i=1}^N a_i(\mu)x^{-2n_i-2} + b(\mu)x^{-2\beta(\mu)-2}\text{ at }\mu_0.
\]
Moreover, on account of Lemma~\ref{lem:lemes_varis} and Definition~\ref{def:moment-finit},
\[
N_n\bigl[[D_{\mu}]^m\bigr]=N_n\bigl[\tfrac{1}{1-z^2}[\CC[D_{\mu}]]_m(\phi^{-1}(z))\bigr]=M_n\bigl[ \CC\bigl[\tfrac{1}{1-z^2}[\CC[D_{\mu}]]_m(\phi^{-1}(z))  \bigr]  \bigr]=M_n\bigl[ [\CC[D_{\mu}]]_m  \bigr].
\]
That is,
\[
N_n\bigl[[D_{\mu}]^m\bigr]=M_n\bigl[ [(\LL_{\boldsymbol\nu_n(\mu)}\circ\CC)[f_{\mu}]]_m  \bigr].
\]
Therefore, if $D_{\mu}$ satisfies one of the hypothesis of the statement, the function $L_{\mu}\!:=(\LL_{\boldsymbol\nu_n(\mu)}\circ\CC)[f_{\mu}]$ satisfies one of the hypothesis in Theorem~\ref{thm:tecnic}. It turns out then that either
\[
(\LL_{\boldsymbol\nu_n(\mu)}\circ\F\circ\CC)[f_{\mu}](x)\sim_{+\infty} C(\mu) x^{\xi(\mu)} \text{ at } \mu_0
\]
or 
\[
(\LL_{\boldsymbol\nu_n(\mu)}\circ\F\circ\CC)[f_{\mu}](x)\sim_{+\infty} C(\mu) \Omega(x,\alpha(\mu))x^{\xi(\mu)} \text{ at } \mu_0
\]
for some functions $C,\xi,\alpha$ with $C(\mu_0)\neq 0$ and $\alpha(\mu_0)=-1$. Let us note that
\begin{align*}
(\LL_{\boldsymbol\nu_n(\mu)}\circ\F\circ\CC)[f_{\mu}](x)&=
\LL_{\boldsymbol\nu_n(\mu)}\bigl[\sqrt{1+x^2}(\CC\circ\F)[f_{\mu}](x) \bigr]\\
&= (\LL_{\boldsymbol\nu_n(\mu)}\circ\CC)\bigl[ (1-z^2)^{-\frac{1}{2}}\F[f_{\mu}](z)\bigr](x)\\
&=(\CC\circ\DD_{\boldsymbol\nu_n(\mu)})\bigl[ (1-z^2)^{-\frac{1}{2}}\F[f_{\mu}](z)\bigr](x),
\end{align*}
with $z=\phi(x)=\frac{x}{\sqrt{1+x^2}}$, where we use $(c)$ in Lemma~\ref{lem:lemes_varis} in the first quality, the identity $\sqrt{1+x^2}\CC[\varphi](x)=\CC[(1-z^2)^{-\frac{1}{2}}\varphi(z)]$ with $\varphi=\F[f_{\mu}]$ in the second equality, and $(b)$ in Lemma~\ref{lem:lemes_varis} in the third equality. So we have
\[
(\CC\circ\DD_{\boldsymbol\nu_n(\mu)})\bigl[ (1-z^2)^{-\frac{1}{2}}\F[f_{\mu}](z)\bigr](x)\sim_{+\infty} C(\mu) \Omega(x,\alpha(\mu))^m x^{\xi(\mu)} \text{ at } \mu_0
\]
with $m\in\{0,1\}$. Using the definition of $\CC$, the previous equation yields to
\[
\frac{1}{1+x^2}\DD_{\boldsymbol\nu_n(\mu)}\bigl[ (1-\phi(x)^2)^{-\frac{1}{2}}\F[f_{\mu}](\phi(x))\bigr]\sim_{+\infty} C(\mu) \Omega(x,\alpha(\mu))^m x^{\xi(\mu)} \text{ at } \mu_0.
\]
Setting $x=\phi^{-1}(z)=\frac{z}{\sqrt{1-z^2}}$, the previous identity implies that
\[
\lim_{(z,\mu)\rightarrow(1,\mu_0)}\frac{(1-z^2)^{1+\frac{\xi(\mu)}{2}}\DD_{\boldsymbol\nu_n(\mu)}\bigl[ (1-z^2)^{-\frac{1}{2}}\F[f_{\mu}](z)\bigr]}{\Omega\bigl(\tfrac{z}{\sqrt{1-z^2}},\alpha(\mu)\bigr)^m}=C(\mu_0).
\]
Thus, on account of the definition of $\DD_{\boldsymbol\nu_n(\mu)}$ in~\eqref{eq:D},
\[
\lim_{(z,\mu)\rightarrow(1,\mu_0)}\frac{(1-z^2)^{1+\frac{\xi(\mu)}{2}}(z(1-z^2))^{\frac{n(n+1)}{2}}}{\Omega\bigl(\tfrac{z}{\sqrt{1-z^2}},\alpha(\mu)\bigr)^m}\frac{W\bigl[\psi_{\nu_1(\mu)}(z),\dots,\psi_{\nu_n(\mu)}(z), (1-z^2)^{-\frac{1}{2}}\F[f_{\mu}](z)\bigr]}{\prod_{i=1}^n \psi_{\nu_i(\mu)}(z)}=C(\mu_0)
\]
which, since $\psi_{\nu}(z)=\frac{z^{\nu}}{(1-z^2)^{1+\nu/2}}$, implies that
\[
\lim_{(z,\mu)\rightarrow(1,\mu_0)}\frac{(1-z^2)^{\kappa(\mu)}}{\Omega\bigl(\tfrac{z}{\sqrt{1-z^2}},\alpha(\mu)\bigr)^m}W\bigl[\psi_{\nu_1(\mu)}(z),\dots,\psi_{\nu_n(\mu)}(z), (1-z^2)^{-\frac{1}{2}}\F[f_{\mu}](z)\bigr]=C(\mu_0)
\]
with $\kappa(\mu)\!:=1+\frac{\xi(\mu)}{2}+\frac{n(n+3)}{2}+\frac{1}{2}\sum_{i=1}^{n}\nu_i(\mu)$ and $m\in\{0,1\}$. The result follows then by Lemma~\ref{lema:crit-finit} and taking the identity~\eqref{FT-finit} into account.
\end{prova}

\section{Applications}\label{sec:applications}

In this Section we present the results obtained applying Theorems~\ref{thm:criticalitat-infinite} and~\ref{thm:criticalitat-finite} to the two-parametric families of centers introduced in Section~\ref{sec:intro}.

\subsection{The family $\ddot{x}+(x+1)^p-(x+1)^q=0$}

The family of vector fields~\eqref{familia-potencial} is an analytic potential system with potential function
\begin{equation}\label{eq:potencial-familia}
V_{\mu}(x)=\int_1^{x+1} (u^p-u^q)du.
\end{equation}
According with the dichotomy that produces $h_0$ in hypothesis~\textnormal{\textbf{(H)}}, in order to prove Theorem~\ref{thm:familia} we split the parameter space in two parts: $\Lambda_1\!:=\{(q,p)\in\Lambda: q>-1\}$ and $\Lambda_2\!:=\{(q,p)\in\Lambda : q<-1\}$. As mentioned in the Introduction, the line $\{q+1=0\}$ correspond to parameters such that in any neighbourhood of them $h_0=+\infty$ and $h_0<+\infty$ coexist. This scenario is out of reach with the current techniques. If $\mu\in\Lambda_1$ system~\eqref{familia-potencial} has a non-degenerate center at the origin and the projection of the period annulus on the $x$-axis is $\PI_{\mu}=(-1,\rho(\mu)-1)$ with $\rho(\mu)\!:=\bigl(\tfrac{p+1}{q+1}\bigr)^{\frac{1}{p-q}}$. Moreover, the energy at the outer boundary is $h_0(\mu)=\tfrac{p-q}{(p+1)(q+1)}<+\infty$ for all $\mu\in\Lambda_1$. On the other hand, if $\mu\in\Lambda_2$ the origin is also a non-degenerate center but in this case $\PI_{\mu}=(-1,+\infty)$ and $h_0(\mu)=+\infty$. (For more details in this direction we refer to~\cite{ManRojVil2017}.) In particular, in both cases hypothesis~\textnormal{\textbf{(H)}} in Definition~\ref{def:Ak} is fulfilled. The proof of Theorem~\ref{thm:familia} follows directly from Propositions~\ref{prop:fam1} and~\ref{prop:resultat-familia}.

\subsubsection{The criticality is at most one in $\Lambda_2$}\label{sec:lambda2}

The purpose of this Section is to apply Theorem~\ref{thm:criticalitat-infinite} with $n=1$ on family~\eqref{familia-potencial} to parameters $\mu_0=(q_0,1)$ with $q_0\in(-3,-1)$, $q_0\neq -2$. According with the Theorem, we first need to compute the asymptotic expression of the function $\LL_{\nu(\mu)}[x(g_{\mu}^{-1})''(x)-x(g_{\mu}^{-1})''(-x)]$ at $x=+\infty$. To this end, we follow the same strategy used on the proof of~\cite[Theorem~A]{ManRojVil2016a}. Since the function $x(g_{\mu}^{-1})''(x)$ is analytic on $\R$, we can write
\begin{equation}\label{eq:suma}
\LL_{\nu(\mu)}[x(g_{\mu}^{-1})''(x)-x(g_{\mu}^{-1})''(-x)]=\frac{W[x^{\nu(\mu)},x(g_{\mu}^{-1})''(x)]}{x^{\nu(\mu)-1}}+\frac{W[x^{\nu(\mu)},-x(g_{\mu}^{-1})''(-x)]}{x^{\nu(\mu)-1}}
\end{equation}
for all $x\in(0,+\infty)$. On account of the equalities $(g^{-1})''(x)=2\mathscr R(g^{-1}(x))$ with $\mathscr R\!:=\frac{(V')^2-2V V''}{(V')^3}$ and $V(g^{-1}(x))=x^2$, equality~\eqref{eq:suma} reads
\begin{equation}\label{eq:suma-S}
\LL_{\nu(\mu)}[x(g_{\mu}^{-1})''(x)-x(g_{\mu}^{-1})''(-x)]= 4 \bigl ( S_{\mu}(g_{\mu}^{-1}(x))+ S_{\mu}(g_{\mu}^{-1}(-x)) \bigr),
\end{equation}
where
\[
S_{\mu}(x)\!:= V_{\mu}'(x)^{-1}V_{\mu}(x)^{2-\frac{\nu(\mu)}{2}}W\bigl[V_{\mu}^{\frac{\nu(\mu)-1}{2}},\mathscr R_{\mu}\bigr](x).
\]

\begin{lema}\label{lema:quant-S}
Let $S_{\mu}$ be defined as above with $V_{\mu}$ defined in~\eqref{eq:potencial-familia}. Taking $\mu_0=(q_0,1)$, the following hold:
\begin{enumerate}[$(a)$]
\item $S_{\mu}(g_{\mu}^{-1}(x))\sim_{+\infty} \frac{p(1+3p)(p-q)}{q+1}(p+1)^{-\frac{3+4p}{p+1}}x^{-\frac{1+3p}{p+1}}+\frac{(3p-2q-1)(1-p+q)(p-q)}{q+1}(p+1)^{\frac{q-4p-2}{p+1}}x^{\frac{1-3p+2q}{p+1}}$ at $\mu_0$.
\item $S_{\mu}(g_{\mu}^{-1}(-x))\sim_{+\infty} \frac{(q-1)(p-q)}{p+1}(-q-1)^{-\frac{2+3q}{q+1}}x^{\frac{1-q}{q+1}}$ at $\mu_0$.
\end{enumerate}
\end{lema}
\begin{prova}
By means of some algebraic manipulations the function $S_{\mu}$ writes
\[
S_{\mu}(x)=\frac{\sqrt{V_{\mu}(x)}}{2V_{\mu}'(x)^5}\psi_{\mu}(x),
\]
where $\psi\!:=-(V'^2-2VV'')((\nu-1)V'^2+6VV'')-4V^2V'V'''$. Using the expression in~\eqref{eq:potencial-familia} we can assert that $\psi_{\mu}$ is the sum of $12$ monomials of the form $c(\mu)(x+1)^{n_1p+n_2q+n_3}$ with $n_i\in\Z$ for $i=1,2,3$ and $c$ a well defined rational function at $\mu=\mu_0$. The biggest exponent for $\mu\approx\mu_0$ is $(x+1)^{4p}$ with coefficient $\frac{(p-1)(p-1+\nu(p+1))}{(p+1)^2}$. Let us notice that, since the coefficient vanishes at $p=p_0=1$, $\psi_{\mu}$ is not continuously quantifiable at infinity in $\mu=\mu_0$ unless we fix $\nu=\frac{p-1}{p+1}$. With this choice of $\nu$, the following two largest exponents for $\mu\approx\mu_0$ are $(x+1)^{3p-1}$ and $(x+1)^{3p+q}$. More precisely, we have
\[
S_{\mu}(x)\sim_{+\infty} \frac{p(1+3p)(p-q)}{(1+p)^{\frac{5}{2}}(1+q)}(x+1)^{-\frac{1}{2}(1+3p)}+\frac{(3p-2q-1)(1+q-p)(p-q)}{(1+p)^{\frac{5}{2}}(1+q)}(x+1)^{-\frac{1}{2}(3p-2q-1)}.
\]
The result in $(a)$ follows using that $g_{\mu}^{-1}(x)\rightarrow +\infty$ as $x\rightarrow +\infty$ and $g_{\mu}(x)^2=V_{\mu}(x)\sim_{+\infty} \frac{(x+1)^{p+1}}{p+1}$. The assertion $(b)$ follows similarly now looking for the smallest exponent for $\mu\approx\mu_0$.
\end{prova}

\begin{prop}\label{prop:fam1}
Let $\{X_{\mu}\}_{\mu\in\Lambda}$ be the family of analytic potential systems~\eqref{familia-potencial} and consider the period function of the center at the origin. If $\mu_0=(q_0,1)$ with $q_0\in(-3,-1)\setminus\{-2\}$ then $\mathrm{Crit}\bigl((\Pi_{\mu_0},X_{\mu_0}),X_{\mu}\bigr)=1$.
\end{prop}
\begin{prova}
From~\cite{ManRojVil2016} we already know that $\mathrm{Crit}\bigl((\Pi_{\mu_0},X_{\mu_0}),X_{\mu}\bigr)\geq 1$. Let us prove the opposite inequality. 
According with equality~\eqref{eq:suma-S} and Lemma~\ref{lema:quant-S}, the function $\LL_{\nu(\mu)}[x(g_{\mu}^{-1})''(x)-x(g_{\mu}^{-1})''(-x)]$ satisfies 
\[
\LL_{\nu(\mu)}[x(g_{\mu}^{-1})''(x)-x(g_{\mu}^{-1})''(-x)]\sim_{+\infty} \frac{p(1+3p)(p-q)}{(p+1)^{\frac{5}{2}}(q+1)}(p+1)^{-\frac{1+3p}{2(p+1)}}x^{-\frac{1+3p}{p+1}}+b(\mu)x^{\beta(\mu)}
\]
at $\mu_0=(q_0,1)$, where
\[
\beta(\mu)\!:=\max_{\mu\approx\mu_0}\left\{\frac{1-3p+2q}{p+1}, \frac{1-q}{q+1} \right\}=\begin{cases}
\frac{1-q}{q+1} & \text{ if } q\in(-3,-2),\\
\frac{1-3p+2q}{p+1} & \text{ if } q\in(-2,-1),\\
\end{cases}
\]
and
\[
b(\mu)\!:=\begin{cases}
\frac{(q-1)(p-q)}{p+1}(-q-1)^{-\frac{2+3q}{q+1}} & \text{ if } q\in(-3,-2),\\
\frac{(3p-2q-1)(1-p+q)(p-q)}{q+1}(p+1)^{\frac{q-4p-2}{p+1}} & \text{ if } q\in(-2,-1).\\
\end{cases}
\]
We point out that $\beta(\mu)$ is continuous but $b(\mu)$ changes sign at $q=-2$. On account of the previous computations, the quantifier at infinity of $\LL_{\nu(\mu)}[x(g_{\mu}^{-1})''(x)-x(g_{\mu}^{-1})''(-x)]$ at $\mu=\mu_0$ is $-2$. In addition, the first momentum $M[[L_{\mu}]_0]$ on Theorem~\ref{thm:criticalitat-infinite} vanishes identically. Indeed,
\[
M[[L_{\mu}]_0]=\int_{-\infty}^{\infty}x(g_{\mu}^{-1})''(x)dx=\int_{-1}^{+\infty}\left(\frac{1}{2}-\frac{V V''}{(V')^2}(x)\right)dx = \left. \frac{V(x)}{V'(x)}-\frac{x}{2}\right|_{-1}^{+\infty}=0
\]
for all $\mu\in\Lambda_2$. Moreover, $\beta(\mu)\in(-3,-2)$ and $b(\mu)$ is continuous for $\mu\approx\mu_0$. Therefore, condition $(b)$ in Theorem~\ref{thm:criticalitat-infinite} holds. This proves that the criticality at $\mu_0$ is at most one.
\end{prova}

\subsubsection{The criticality is at most one in $\Lambda_1$}\label{sec:lambda1}

The energy at the outer boundary of the period annulus for parameters $\mu\in\Lambda_1$ is finite. Accordingly we shall use Theorem~\ref{thm:criticalitat-finite} with $n=1$ to parameters $\mu_0=(q_0,p_0)$ with $p_0+2q_0+1=0$, $q_0\in(-1,-\tfrac{1}{2})$ in order to prove the second assertion of $(b)$ in Theorem~\ref{thm:familia}. A similar argument as in the previous section shows that
\begin{equation}\label{eq:suma2}
\DD_{\nu(\mu)}[z\sqrt{h_0}(g_{\mu}^{-1})''(z\sqrt{h_0})-z\sqrt{h_0}(g_{\mu}^{-1})''(-z\sqrt{h_0})]=\frac{4}{h_0}\left(S_{\mu}\bigl(g_{\mu}^{-1}(z\sqrt{h_0})\bigr)-S_{\mu}\bigl(g_{\mu}^{-1}(-z\sqrt{h_0})\bigr) \right),
\end{equation}
where
\[
S_{\mu}(x)\!:=\frac{W\left[\bigl(\frac{V_{\mu}}{h_0-V_{\mu}}\bigr)^{\frac{\nu}{2}}, (h_0-V_{\mu})V_{\mu}^{\frac{1}{2}}\mathscr R_{\mu}\right](x)}{(h_0-V_{\mu}(x))^{-\frac{\nu}{2}}V_{\mu}(x)^{\frac{\nu}{2}-1}V_{\mu}'(x)}.
\]
Here, and in what follows, we omit the dependence in $\mu$ of $h_0$ for the sake of simplicity. Next result is general for any potential function $V_{\mu}$ with finite energy at the outer boundary.

\begin{lema}\label{lema:quant-S-dreta}
Let $S_{\mu}$ be defined as above and fix $\mu_0\in\Lambda$. Assume that, for all $\mu$ in a neighbourhood of $\mu_0$, $h_0(\mu)$ is finite and the right endpoint of the projection of the period annulus $x_r=x_r(\mu)$ is also finite and satisfies $V_{\mu}'(x_r)\neq 0$. Then
\[
S_{\mu}\bigl( g_{\mu}^{-1}(z\sqrt{h_0})\bigr)\sim_1 \dfrac{h_0^{\frac{3}{2}}(2+\nu)(2h_0V_{\mu}''(x_r)-V_{\mu}'(x_r)^2)}{2V_{\mu}'(x_r)^3}+c(\mu)(1-z^2) \text{ at }\mu_0
\]
with 
\[
c(\mu)\!:=\frac{h_0^{\frac{3}{2}}(12h_0^2(4+\nu)V_{\mu}''(x_r)^2-4h_0^2(4+\nu)V_{\mu}'''(x_r)V_{\mu}'(x_r)-8h_0(5+\nu)V_{\mu}''(x_r)V_{\mu}'(x_r)^2+(8+\nu)V_{\mu}'(x_r)^4)}{4V_{\mu}'(x_r)^5}.
\]
The analogous result is true for $x_{\ell}$.
\end{lema}

\begin{prova}
The function $S_{\mu}$ can be written, after some algebraic manipulations, as
\begin{equation}\label{S-simple}
S_{\mu}(x)=\frac{\sqrt{V_{\mu}(x)}}{2V_{\mu}'(x)^5}\psi_{\mu}(x),
\end{equation}
where $\psi\!:=-\bigl((V')^2-2VV''\bigr)\bigl((V')^2(h_0(\nu-1)+3V)+6(h_0-V)VV''\bigr)-4V^2(h_0-V)V'V'''$. Since the function $V_{\mu}(x)$ is analytic at $x=x_r$ the result follows by considering the Taylor's expansion at $x=x_r$ of the previous expression and using the change of variable $z=g_{\mu}(x)/\sqrt{h_0}$.
\end{prova}

\begin{lema}\label{lema:quant-S-finit}
Let $S_{\mu}$ be defined as above with $V_{\mu}$ defined in~\eqref{eq:potencial-familia}. Taking $\nu(\mu)\equiv -1$ and $\mu_0=(q_0,p_0)$ with $p_0+2q_0+1=0$ and $q_0\in(-1,-\frac{1}{2})$ then
\[
S_{\mu}\bigl( g_{\mu}^{-1}(-z\sqrt{h_0})\bigr)\sim_1  -\dfrac{h_0^{-\frac{1+3q}{2(1+q)}}(p-q)^2q(1+q)^{-\frac{4+5q}{1+q}}(1+3q)}{(1+p)^2}(1-z^2)^{-\frac{1+2q}{1+q}} \text{ at }\mu_0.
\]
Moreover, the asymptotic expression in Lemma~\ref{lema:quant-S-dreta} holds with $c(\mu_0)\neq 0$.
\end{lema}
\begin{prova}
The first assertion of the result can be proven similarly as Lemma~\ref{lema:quant-S} so we skip the details for the sake or brevity. We only prove that $c(\mu_0)\neq 0$. Indeed, using the expression in~\eqref{eq:potencial-familia} and substituting $p_0=-2q_0-1$ some tedious but elementary computations show that (see Lemma~\ref{lema:quant-S-dreta})
\begin{align*}
c(\mu_0)&=\frac{16^{-\frac{q}{1+3q}}h_0^{\frac{3}{2}}(1+3q)^2\bigl(-1+\tfrac{1}{1+q}\bigr)^{-\frac{4q}{1+3q}}}{4q(1+q)^4V_{\mu_0}'(\rho)^5}\left(-7q(1+3q)^2+4^{\frac{3+8q}{1+3q}}q(1+q)(1+3q)\bigl(-1+\tfrac{1}{1+q}\bigr)^{\frac{1+q}{1+3q}}\right. \\
&\phantom{=}\left.+9\times 2^{\frac{4+8q}{1+3q}}(1+q)^2\bigl(-1+\tfrac{1}{1+q}\bigr)^{\frac{2(1+q)}{1+3q}} \right).
\end{align*}
We stress that each term inside the parenthesis is positive for $q_0\in(-1,-\tfrac{1}{2})$. Also the multiplicative term is non-vanishing so $c(\mu_0)\neq 0$. 
\end{prova}

The following result is the version of~\cite[Lemma~3.12]{ManRojVil2016a} for $\ell=n=1$.

\begin{lema}\label{lema:redu}
Let $f\in C^{\omega}[0,1)$ and $\nu\in\R$. Let us assume that $f$ is quantifiable at $z=1$ by $\xi$. If $\xi<\frac{1}{2}$ then
\[
N\bigl[\DD_{\nu}[f]\bigr]=-(1+\nu)N[f].
\]
\end{lema}

\begin{prop}\label{prop:resultat-familia}
Let $\{X_{\mu}\}_{\mu\in\Lambda}$ be the family of analytic potential systems~\eqref{familia-potencial} and consider the period function of the center at the origin. If $\mu_0=(q_0,p_0)$ with $p_0+2q_0+1=0$ and $q_0\in(-1,-\tfrac{1}{2})\setminus\{-\tfrac{2}{3}\}$ then $\mathrm{Crit}\bigl((\Pi_{\mu_0},X_{\mu_0}),X_{\mu}\bigr)=1$.
\end{prop}

\begin{prova}
From~\cite{ManRojVil2016} we already know that $\mathrm{Crit}\bigl((\Pi_{\mu_0},X_{\mu_0}),X_{\mu}\bigr)\geq 1$. Let us prove the opposite inequality. 
Let us denote $f_{\mu}(z)\!:=z\sqrt{h_0}(g_{\mu}^{-1})''(z\sqrt{h_0})-z\sqrt{h_0}(g_{\mu}^{-1})''(-z\sqrt{h_0})$ for the sake of simplicity. In \cite[Lemma~4.4]{ManRojVil2016a} it was shown that $f_{\mu}$ satisfies the hypothesis on Lemma~\ref{lema:redu} for all $\mu\approx\mu_0$. Applying Lemma~\ref{lema:redu} with $\nu(\mu)\equiv -1$ we have then $N\bigl[\DD_{\nu(\mu)}[f_{\mu}]\bigr]=0$ for all $\mu\approx\mu_0$. Moreover, on account of the equality~\eqref{eq:suma2} and Lemma~\ref{lema:quant-S-finit},
\[
\DD_{\nu(\mu)}[f_{\mu}](z)\sim_{1} \dfrac{h_0^{\frac{3}{2}}(2h_0V_{\mu}''(\rho)-V_{\mu}'(\rho)^2)}{2V_{\mu}'(\rho)^3}+ b(\mu)(1-z^2)^{\beta(\mu)} \text{ at }\mu_0,
\]
where
\[
\beta(\mu)\!:=\min_{\mu\approx\mu_0}\left\{1,-\frac{1+2q}{1+q} \right\} = 
\begin{cases}
1 & \text{ if } q\in(-1,-\tfrac{2}{3}),\\
-\frac{1+2q}{1+q} & \text{ if } q\in(-\tfrac{2}{3},-\tfrac{1}{2}),
\end{cases}
\]
and
\[
b(\mu)\!:=\begin{cases}
c(\mu) & \text{ if } q\in(-1,-\tfrac{2}{3}),\\
-\dfrac{h_0^{-\frac{1+3q}{2(1+q)}}(p-q)^2q(1+q)^{-\frac{4+5q}{1+q}}(1+3q)}{(1+p)^2} & \text{ if } q\in(-\tfrac{2}{3},-\tfrac{1}{2}),
\end{cases}
\]
due to $c(\mu_0)\neq 0$ by Lemma~\ref{lema:quant-S-finit}. Let us apply now Theorem~\ref{thm:criticalitat-finite} with $n=1$. According with the previous discussion, if $p_0+2q_0+1=0$ and $q_0\in(-\tfrac{2}{3},-\tfrac{1}{2})$ then assertion in $(b)$ of Theorem~\ref{thm:criticalitat-finite} is satisfied with $N=1$, $n_1=0$ and $\ell=1$ due to $N\bigl[\DD_{\nu(\mu)}[f_{\mu}(z)]\bigr]\equiv 0$ and $\beta(\mu)-1\leq\tfrac{1}{2}$. Then $\mathrm{Crit}\bigl((\Pi_{\mu_0},X_{\mu_0}),X_{\mu}\bigr)\leq 1$ in this case. On the other hand, if $q_0\in(-1,-\tfrac{2}{3})$ we apply Theorem~\ref{thm:criticalitat-finite} with $N=2$, $n_1=0$ and $n_2=1$. By Lemma~\ref{lema:moment2} we also have
\[
N\left[\bigl[\DD_{\nu}[f_{\mu_0}]\bigr]^1 \right] \neq 0.
\]
Then assertion $(a)$ in Theorem~\ref{thm:criticalitat-finite} states that $\mathrm{Crit}\bigl((\Pi_{\mu_0},X_{\mu_0}),X_{\mu}\bigr)\leq 1$ also in this case.
\end{prova}

\subsection{The family of dehomogenized Loud's centers}

\begin{figure}[t]
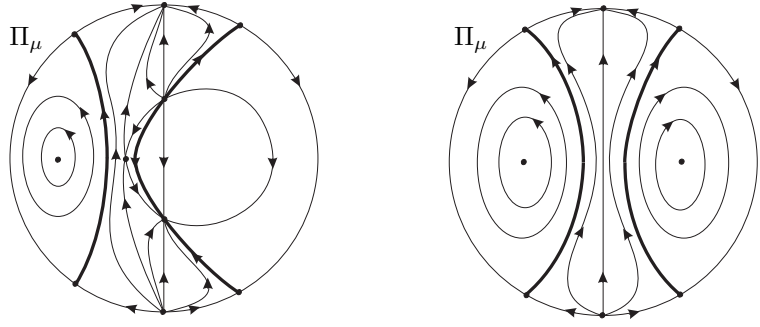

\centering
\begin{lpic}[l(0mm),r(0mm),t(0mm),b(5mm)]{dib3(2)}
\lbl[l]{0,18.5;$\Pi_{\mu}$}
\lbl[l]{29.25,18.5;$\Pi_{\mu}$}
\end{lpic}
\caption{\label{fig:fp-loud} Phase portrait of~\eqref{loud} in the Poincaré disc for $\mu=(D,F)\in\Lambda$ with $D<-1$ (left) and $D>-1$ (right). The center at $(0,0)$ is placed on the left of the centered invariant line $\{x=1\}$ for convenience. The invariant hyperbola $\mathscr C_{\mu}$ is in boldface type. (Figure extracted from~\cite{MarMarVil2006}.)}
\end{figure}

Through all this section we consider parameters $\mu=(D,F)$ inside the open set
\[
\Lambda\!:=\{(D,F)\in\R^2 : 1<F<\tfrac{3}{2}, D<-\tfrac{1}{2}, D+F>0\}.
\]
The Loud system~\eqref{loud} has a first integral given by
\[
H_{\mu}(x,y)=(1-x)^{-2F}\left(\frac{1}{2}y^2-q_{\mu}(x)\right)
\]
for all $F\notin\{0,1,\tfrac{1}{2}\}$ (see for instance~\cite{MarMarVil2006}) where $q_{\mu}(x)=a(\mu)x^2+b(\mu)x+c(\mu)$ with
\[
a=\frac{D}{2(1-F)},\ b=\frac{D-F+1}{(1-F)(1-2F)} \text{ and }c=\frac{F-D-1}{2F(1-F)(1-2F)},
\]
and integrating factor $\kappa(x)=(1-x)^{-2F-1}$. The line at infinity $L_{\infty}$, the line $\{x=1\}$ and the conic $\mathscr C_{\mu}=\{\tfrac{1}{2}y^2-q_{\mu}(x)=0\}$ are invariant curves of the differential system and, for parameters $\mu\in\Lambda$, $\mathscr C_{\mu}$ is a hyperbola intersecting the $x$-axis at
\begin{equation}\label{def:p1p2}
x=p_1(\mu)\!:=\frac{-b-\sqrt{b^2-4ac}}{2a} \text{ and } x=p_2(\mu)\!:=\frac{-b+\sqrt{b^2-4ac}}{2a},
\end{equation}
with $0<p_1(\mu)<p_2(\mu)$. The outer boundary of the period annulus of the center at the origin of system~\eqref{loud} consists of the branch of the hyperbola $\mathscr C_{\mu}$ passing through $(p_1,0)$ and the line at infinity $L_{\mu}$, joined by two hyperbolic saddles (see Figure~\ref{fig:fp-loud}.) Although it is not relevant at this moment (but it will be nearly soon) we point out that there is a bifurcation on the phase portrait at $D=-1$ due to the fact that the branch of hyperbola passing through $(p_2,0)$ crosses the invariant line $\{x=1\}$. In particular, $p_2(\mu)>1$ if and only if $D>-1$.

In~\cite{MarMarVil2006} the authors give an implicit expression of the bifurcation curve $D=\mathcal G(F)$ in $\Lambda$. More concretely, for each $\mu\in\Lambda$ and $s>0$ small let $P(s;\mu)$ be the period of the periodic orbit of system~\eqref{loud} passing through the point $(p_1(\mu)-s,0)$. Then, for parameters in $\Lambda$, the derivative of the period function $P'(s;\mu)$ tends to $\Delta(\mu)$ as $s$ tends to zero uniformly on compact subsets of $\Lambda$, where
\[
\Delta(\mu)\!:=\frac{-1}{\sqrt{2a}(p_2-p_1)(1-p_1)}\left\{ 2-\int_0^1 (1-u)^{-\frac{3}{2}}\left( u^{2-2F} \left(\frac{1-p_2}{1-p_1}(u-1)+1\right)^{2F-1}-1\right)du \right\}.
\]
The bifurcation curve $D=\mathcal G(F)$ is given by the zero set of $\Delta(\mu)$ and the following properties are deduced:
\begin{enumerate}[$(a)$]
\item $-F<\mathcal G(F)<-\tfrac{1}{2}$ for all $F\in(1,\tfrac{3}{2})$.
\item $\lim_{F\rightarrow \tfrac{3}{2}} \mathcal G(F)=-\tfrac{3}{2}$, and
\item $\lim_{F\rightarrow 1} \mathcal G(F)=-\tfrac{1}{2}$.
\end{enumerate}
For the proof of these results we refer to \cite[Theorem~3.6 and Proposition~3.11]{MarMarVil2006}. In particular the curve $D=\mathcal G(F)$ is analytic on $\Lambda$ and joins the parameters $\mu=(-\tfrac{3}{2},\tfrac{3}{2})$ and $\mu=(-\tfrac{1}{2},1)$. The following result gives a more compact expression of the zero set of $\Delta(\mu)$ and an additional property. The proof is given in the Appendix.

\begin{prop}\label{prop:corba-bif-expr}
For parameters $\mu\in\Lambda$ the zero sets of $\Delta(\mu)$ and
$\tfrac{1}{\Gamma\bigl(\tfrac{7}{2}-4F\bigr)}{}_2F_1\left(-\tfrac{1}{2},\tfrac{3}{2};\tfrac{5}{2}-2F;\tfrac{p_2-1}{p_2-p_1}\right)$ coincide. Moreover, $\lim_{F\rightarrow \tfrac{5}{4}} \mathcal G(F)=-1$.
\end{prop}

The Loud system~\eqref{loud} is not a potential system and, at first glance, is not suitable to use the techniques developed in this paper to study the criticality at the outer boundary of its center. However we recall that it has a first integral quadratic in $y$ and its integrating factor depends only on $x$. These two properties are the requirements of~\cite[Lemma~14]{GGV}, which shows that the change of variables
\[
(u,v)=(\phi(1-x),(1-x)^{-F}y), \text{ with }\phi(z)\!:=\frac{z^{-F}-1}{F},
\]
transforms system~\eqref{loud} to the potential system
\begin{equation}\label{loud-potential}
\begin{cases}
\dot u &= -v,\\
\dot v &= (Fu+1)\bigl((Fu+1)^{-\frac{1}{F}}-1\bigr)\bigl(D(Fu+1)^{-\frac{1}{F}}-D-1\bigl).
\end{cases}
\end{equation}
The potential system above has a non-degenerated center at the origin and all the properties that we derive of its period function are directly transmitted to the period function of the Loud's center. In these new variables, the projection of the period annulus is $\mathcal I_{\mu}\!:=(-\tfrac{1}{F},u_r(\mu))$, where
\begin{equation}\label{ur_loud}
u_r(\mu)\!:=\phi(1-p_1(\mu))=\frac{(1-p_1(\mu))^{-F}-1}{F}.
\end{equation}
Let $H_{\mu}(u,v)=\tfrac{1}{2}v^2+V_{\mu}(u)$ by the Hamiltonian function associated to system~\eqref{loud-potential} with $V_{\mu}(0)=0$. If we set $z=\phi^{-1}(u)=(Fu+1)^{-1/F}$ the following expressions will be useful to simplify the forthcoming computations:
\begin{equation}\label{V-loud}
\begin{array}{ll}
V_{\mu}(u)=h_0(\mu)-z^{-2F}V_0(z;\mu) & \text{ with } V_0(z;\mu)=\frac{D}{2-2F}z^2+\frac{1+2D}{2F-1}z-\frac{D+1}{2F},\\
V_{\mu}'(u)=z^{-F}V_1(z;\mu) & \text{ with } V_1(z;\mu)=(z-1)\bigl(D(z-1)-1),\\
V_{\mu}''(u)=V_2(z;\mu) & \text{ with } V_2(z;\mu)=D(F-2)z^2-(2D+1)(F-1)z+F(D+1),\\
V_{\mu}'''(u)=z^F V_3(z;\mu) & \text{ with } V_3(z;\mu)=-2D(F-2)^2+(2D+1)(F-1)z.
\end{array}
\end{equation}
Here $h_0(\mu)=\frac{F-D-1}{2F(F-1)(2F-1)}$ denotes the energy of the outer boundary of the center for all $\mu\in\Lambda$. In particular, $h_0(\mu)$ is finite in $\Lambda$.

\subsubsection{Quantification of the functions involved}

With the aim of applying Theorem~\ref{thm:criticalitat-finite} with $n=1$ to parameters $\mu=(D,F)\in\Lambda$ satisfying $D=\mathcal G(F)$, we need to compute the first terms of the asymptotic development of $D_{\mu}(z)$ at $z=1$. We recall that the previous is reduced to study the asymptotic development of
\[
S_{\mu}(u)\!:=\frac{W\left[\bigl(\frac{V_{\mu}}{h_0-V_{\mu}}\bigr)^{\frac{\nu}{2}}, (h_0-V_{\mu})V_{\mu}^{\frac{1}{2}}\mathscr R_{\mu}\right]\!(u)}{(h_0-V_{\mu}(u))^{-\frac{\nu}{2}}V_{\mu}(u)^{\frac{\nu}{2}-1}V_{\mu}'(u)}
\]
at $u=-\tfrac{1}{F}$ and $u=u_r(\mu)$ due to equality~\eqref{eq:suma2}. Here, and in what follows, we omit the dependence in $\mu$ of $h_0$ and $\nu$ for the sake of brevity. Moreover we use $u$ instead of $x$ to be consistent with the notation in~\eqref{loud-potential}. 

\begin{lema}\label{lema:quantificar_S_loud}
Let $S_{\mu}$ and $V_{\mu}$ be defined as above. Taking $\nu\equiv -1$ we have that
\[
S_{\mu}\bigl( g_{\mu}^{-1}(-z\sqrt{h_0})\bigr)\sim_1 a_1(\mu)(1-z^2)^{-\frac{3}{2}+\frac{1}{2(F-1)}} + a_2(\mu)(1-z^2)^{\frac{2-F}{2F-2}} \text{ at }\mu_0
\]
where
\[
a_1\!:=\frac{D^{\frac{F}{2-2F}}(1+D-F)^2(2-F)(2F-3)h_0^{-1+\frac{1}{2(F-1)}}}{4F^2(2F-1)^2(F-1)^3(2-2F)^{\frac{3}{2}-\frac{1}{2(F-1)}}} \text{ and }
a_2\!:=\frac{D^{\frac{F}{2-2F}}(1+D-F)(3-F)h_0^{\frac{1}{2(F-1)}}}{4F(2F-1)(F-1)^3(2-2F)^{\frac{2-F}{2(F-1)}}}.
\]
\end{lema}
\begin{prova}
In order to get the desired asymptotic expression of $S_{\mu}\bigl( g_{\mu}^{-1}(-z\sqrt{h_0})\bigr)$ at $z=1$ we study the asymptotic expression of $S_{\mu}(u)$ at $u=u_{\ell}=-1/F$. The result will follow then on account of the change of variable $z=-g_{\mu}(u)/\sqrt{h_0}$. To do so we invoke the expression in~\eqref{S-simple} and perform the change of variable $z=\phi^{-1}(u)=(Fu+1)^{-1/F}$. Taking advantage of the relations in~\eqref{V-loud} we can write
$S_{\mu}(\phi(z))$ in terms of the polynomials $V_i(z;\mu)$ and some powers of $z$ depending on the parameter $F$. Due to the fact that $\phi(z)\rightarrow -1/F$ as $z\rightarrow +\infty$, we are concerned about the asymptotic expression of $S_{\mu}(\phi(z))$ at $z=+\infty$. An analogous study than the one done in Lemma~\ref{lema:quant-S} shows in this case that
\[
S_{\mu}(\phi(z))\sim_{+\infty} \frac{(2-F)(D+1-F)^2(2F-3)\sqrt{h_0}}{4F^2D^2(F-1)^3(2F-1)^2}z^{3F-4}+\frac{(3-F)(D+1-F)\sqrt{h_0}}{4FD(F-1)^3(2F-1)}z^{-2+F} \text{ at }\mu_0.
\]
The result follows undoing the change of variable.
\end{prova}

\subsubsection{Proof of Theorem~\ref{thm:loud}}\label{sec:proof_ThmB}

From the results in~\cite{MarMarVil2006} we already know that $\mathrm{Crit}\bigl((\Pi_{\mu_0},X_{\mu_0}),X_{\mu}\bigr)\geq 1$ for parameters $\mu_0$ satisfying $D_0=\mathcal G(F_0)$. The strategy for showing the opposite inequality is to apply Theorem~\ref{thm:criticalitat-finite} with $n=1$ to the family~\eqref{loud-potential}. To do so, let $f_{\mu}$ be defined as in the statement of the Theorem. The first thing to notice is that $f_{\mu}$ satisfies hypothesis~\textnormal{\textbf{(H)}} with $h_0(\mu)<+\infty$ for all $\mu\in\Lambda$. Let us fix $\nu(\mu)\equiv -1$. On account of equality~\eqref{eq:suma2} and Lemmas~\ref{lema:quant-S-dreta} and~\ref{lema:quantificar_S_loud} we have that
\[
\DD_{\nu(\mu)}[f_{\mu}](z)\sim_{1} \dfrac{h_0^{\frac{3}{2}}(2h_0V_{\mu}''(\rho)-V_{\mu}'(\rho)^2)}{2V_{\mu}'(\rho)^3}+ b(\mu)(1-z^2)^{\beta(\mu)} \text{ at }\mu_0,
\]
where
\[
\beta(\mu)\!:=\min_{\mu\approx\mu_0}\left\{1,-\frac{3}{2}+\frac{1}{2(F-1)} \right\} = 
\begin{cases}
1 & \text{ if } F\in(1,\tfrac{6}{5}),\\
-\frac{3}{2}+\frac{1}{2(F-1)} & \text{ if } F\in(\tfrac{6}{5},\tfrac{4}{3}),
\end{cases}
\]
and
\[
b(\mu)\!:=\begin{cases}
c(\mu) & \text{ if } F\in(1,\tfrac{6}{5}),\\
a_1(\mu) & \text{ if } F\in(\tfrac{6}{5},\tfrac{4}{3}).
\end{cases}
\]

If $F_0\in[\tfrac{5}{4},\tfrac{4}{3})$ then $\beta(\mu_0)\leq 1/2$ and assertion $(b)$ in Theorem~\ref{thm:criticalitat-finite} is fulfilled with $N=1$, $n_1=0$ and $\ell=1$. Therefore, assertion $(a)$ in Theorem~\ref{thm:loud} holds.

If $F_0\in(1,\tfrac{5}{4})$ then $\beta(\mu_0)>1/2$. Let us split this interval in two parts, namely $(1,\tfrac{6}{5})$ and $(\tfrac{6}{5},\tfrac{5}{4})$. The parameter $F_0=\tfrac{6}{5}$ is excluded due to the noncontinuity of the function $\beta(\mu)$. If $F_0\in(\tfrac{6}{5},\tfrac{5}{4})$ then the coefficient of $(1-z^2)^{\beta(\mu)}$ in the previous asymptotic expansion is the nonvanishing function $a_1(\mu)$ on Lemma~\ref{lema:quantificar_S_loud}. Moreover, by Lemma~\ref{lema:redu}, $N[\DD_{\nu(\mu)}[f_{\mu}]]=0$ for all $\mu\approx\mu_0$. Consequently,  we apply Theorem~\ref{thm:criticalitat-finite} with $N=2$, $n_1=0$ and $n_2=1$. Using Lemma~\ref{lema:moment2-loud} we have that
\[
N\!\left[\bigl[\DD_{\nu}[f_{\mu_0}]\bigr]^1 \right]=0 \text{ if and only if }\frac{{}_2F_1\bigl(-\frac{3}{2},\frac{5}{2};\frac{7}{2}-4F_0;\alpha(\mu_0)\bigr)}{\Gamma\bigl(\tfrac{7}{2}-4F_0\bigr)}=0.
\]
Proposition~\ref{prop:casitot} ensures that the previous Hypergeometric function do not vanishes for $F_0\in(\tfrac{9}{8},\tfrac{5}{4})$. Consequently, hypothesis of Theorem~\ref{thm:criticalitat-finite} are fulfilled and so the criticality at the outer boundary for those parameters is exactly one.

Let us finally consider $F_0\in(1,\tfrac{6}{5})$. In this case the coefficient of $(1-z^2)^{\beta(\mu)}$ on the asymptotic expansion above is $c(\mu)$ given in Lemma~\ref{lema:quant-S-dreta}. Here then we invoke the assumption in the statement of Theorem~\ref{thm:loud} regarding $c(\mu_0)\neq 0$. As before we can apply Theorem~\ref{thm:criticalitat-finite} with $N=2$, $n_1=0$ and $n_2=1$ and the hypothesis are satisfied whenever condition~\eqref{eq:condition} is fulfilled. That is, whenever the hypergeometric function above do not vanish. Again Proposition~\ref{prop:casitot} ensures that this is true for parameters $F_0\in(\tfrac{9}{8},\tfrac{5}{4})$. Consequently, the result in assertion $(b)$ of Theorem~\eqref{thm:loud} holds. For parameters $F_0\in(1,\tfrac{9}{8}]$ condition~\eqref{eq:condition} is required and so the result follows if it is satisfied. This ends the proof of the first assertion in Theorem~\ref{thm:loud}.

\subsubsection{Some additional comments}\label{sec:comments}

This section is devoted to discuss some technicalities involved in the proof of Theorem~\ref{thm:loud}. More concretely, we show that condition~\eqref{eq:condition} is fulfilled when $F_0\in(\tfrac{9}{8},\tfrac{5}{4})$ and we also give some numerical intuition to the fact that the condition should be fulfilled for every $F_0\in(1,\tfrac{5}{4})$. 

\begin{prop}\label{prop:casitot}
If $\mu=(D,F)$ satisfy $D=\mathcal G(F)$ with $\tfrac{9}{8}<F<\tfrac{5}{4}$ then 
\[
\frac{{}_2F_1\bigl(-\frac{3}{2},\frac{5}{2};\frac{7}{2}-4F;\alpha(\mu)\bigr)}{\Gamma\bigl(\tfrac{7}{2}-4F\bigr)}\neq 0.
\] 
\end{prop}

\begin{prova}
From the fact that the parameters $\mu=(D,F)$ satisfy $D=\mathcal G(F)$ with $\tfrac{9}{8}<F<\tfrac{5}{4}$ it follows that $D>-1$ and so $p_2(\mu)>1$. Moreover $p_1<1$. Consequently, $\alpha(\mu)=\tfrac{p_2-1}{p_2-p_1}\in(0,1)$.

Consider the function
\[
\varphi:=(z,F)\longrightarrow \frac{{}_2F_1\bigl(-\frac{3}{2},\frac{5}{2};\frac{7}{2}-4F;z\bigr)}{\Gamma\bigl(\tfrac{7}{2}-4F\bigr)}.
\]
From~\cite{AS} the derivative of $\varphi(z,F)$ with respect to $z$ writes
\[
\frac{d\varphi}{dz}(z,F)=-\frac{15}{4\Gamma\bigl(\tfrac{9}{2}-4F\bigr)}{}_2F_1\bigl(-\frac{1}{2},\frac{7}{2};\frac{9}{2}-4F;z\bigr).
\]
Using the series expression of the Hypergeometric function it turns out that all its coefficients are positive for $F\in(\tfrac{9}{8},\tfrac{5}{4})$. Consequently, for each $F\in(\tfrac{9}{8},\tfrac{5}{4})$ we have that
\[
\frac{d\varphi}{dz}(z,F)>0 \text{ if }z>0.
\]
Hence for each $F\in(\tfrac{9}{8},\tfrac{5}{4})$ the function $\varphi(\cdot,F)$ is increasing for $z>0$. In addition, 
\[
\varphi(0,F)=\frac{1}{\Gamma\bigl(\tfrac{7}{2}-4F\bigr)}>0
\]
if $F\in(\tfrac{9}{8},\tfrac{5}{4})$. Therefore, for each $F\in(\tfrac{9}{8},\tfrac{5}{4})$, the function $\varphi(z,F)$ is positive for all $z>0$. Taking $z=\alpha(\mu)$ and on account of $0<\alpha(\mu)<1$ for all parameters under consideration, the result follows.
\end{prova}

\begin{figure}
\centering
\subfloat[\label{fig:tall}]
{
 \includegraphics[scale=.5]{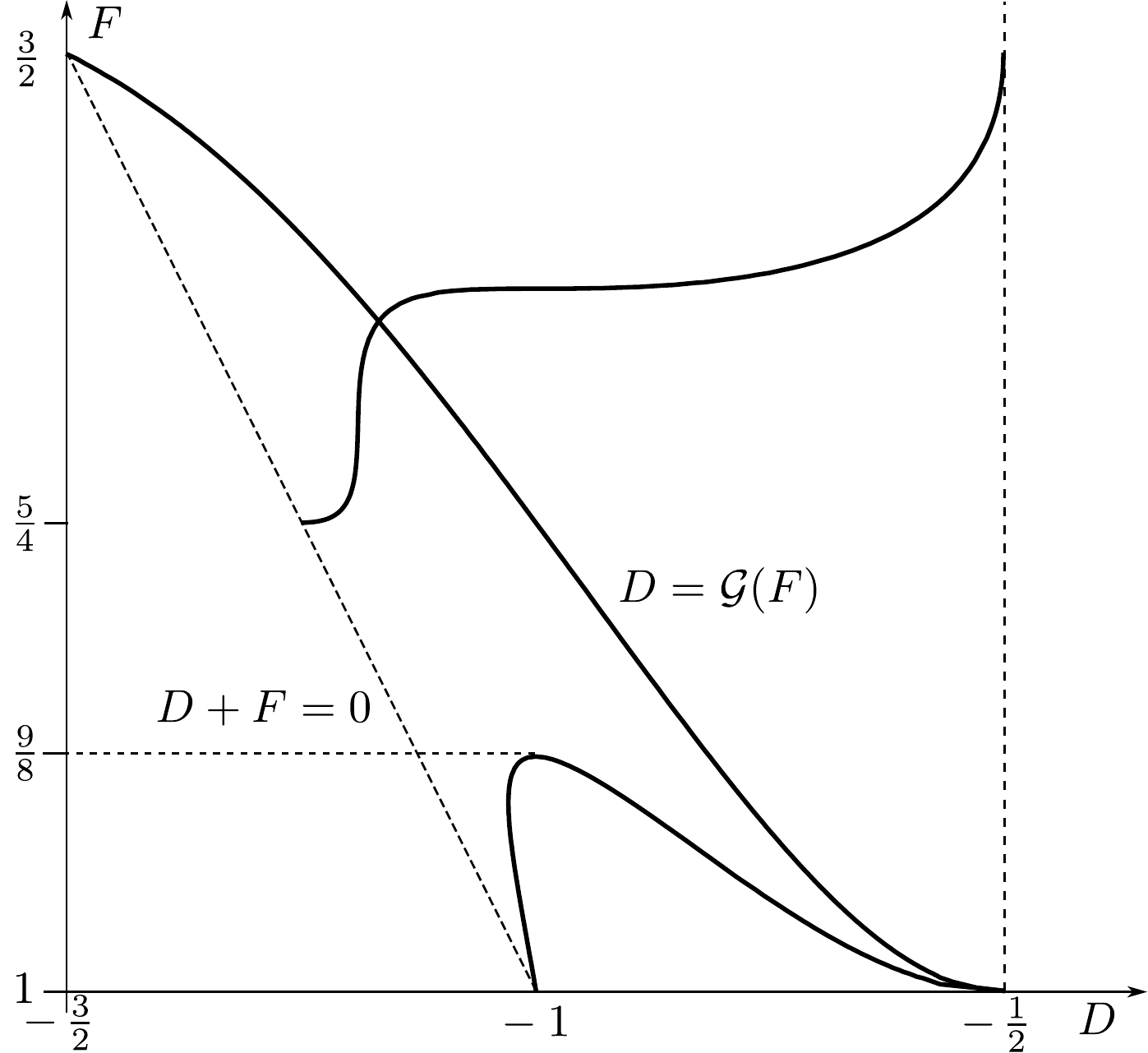}
}
\hspace{1cm}
\subfloat[\label{fig:zoom}]
{
\includegraphics[scale=.5]{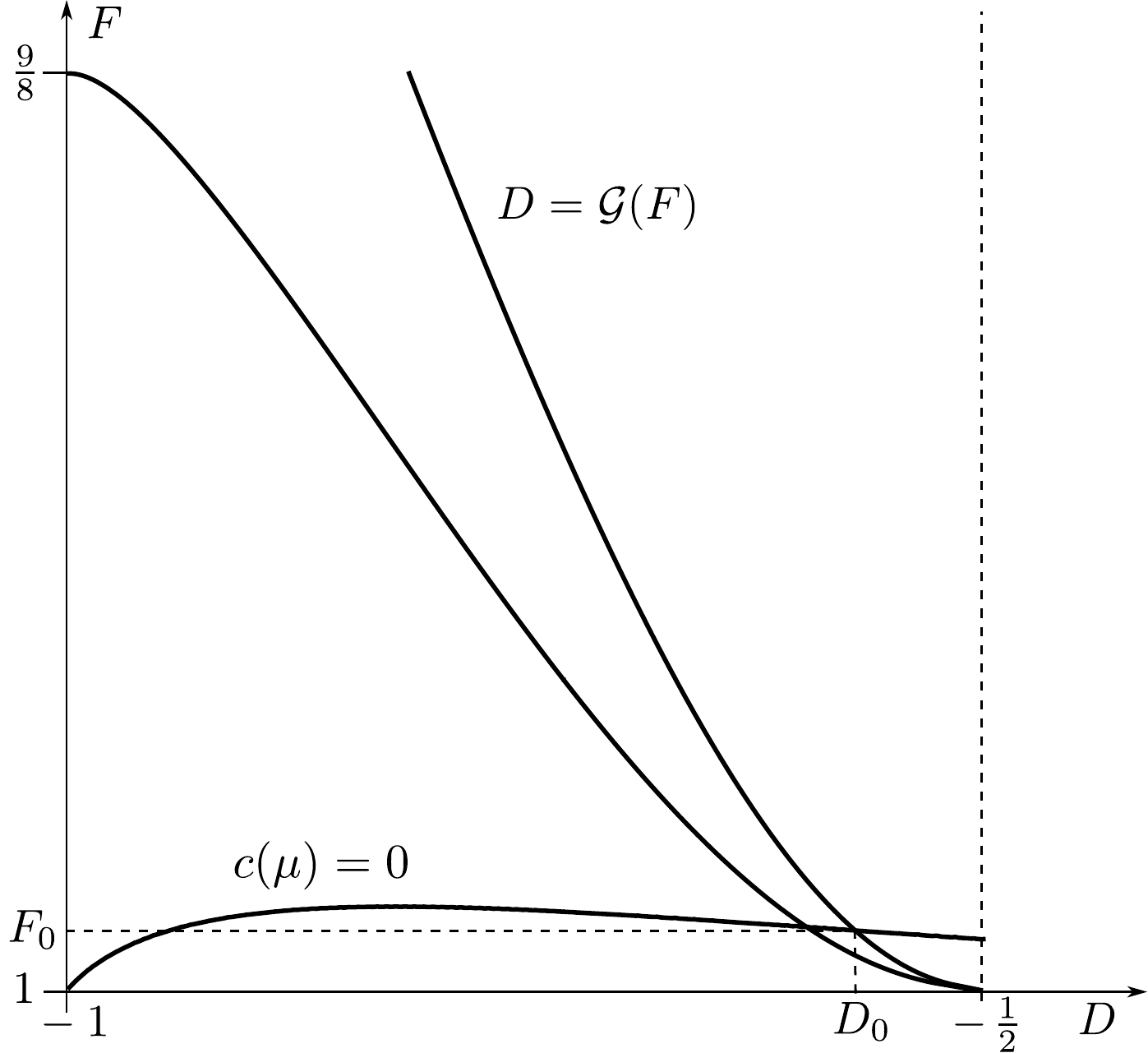}
}
\caption{(a) The curve joining the points $(-\tfrac{3}{2},\tfrac{3}{2})$ and $(-\tfrac{1}{2},1)$ correspond to the bifurcation parameters $D=\mathcal G(F)$. The other two bold lines are the two connected components of the zero set of the function~\eqref{eq:second}.  (b) Zoom of the previous figure in the parameter region $(D,F)\in(-1,-\tfrac{1}{2})\times(1,\tfrac{9}{8})$ and the curve $c(\mu)=0$.}
\end{figure}

Figure~\ref{fig:tall} exhibits the curve $D=\mathcal G(F)$ together with the two connected components of the zero set level of the function 
\begin{equation}\label{eq:second}
\frac{{}_2F_1\bigl(-\frac{3}{2},\frac{5}{2};\frac{7}{2}-4F;\alpha(\mu)\bigr)}{\Gamma\bigl(\tfrac{7}{2}-4F\bigr)}.
\end{equation}
As the previous result states, there is no crossing of the zero set with the curve $D=\mathcal G(F)$ for parameters with $F\in(\tfrac{9}{8},\tfrac{5}{4})$. Indeed Proposition~\ref{prop:casitot} shows that condition~\eqref{eq:condition} is satisfied in this situation. Numerics seems to show that the non-crossing property is still verified till $F=1$. However, we have not been able to prove it analytically.

We end this section with also a numerical intuition about the equation $c(\mu_0)=0$ in the statement of Theorem~\ref{thm:loud}. Although again no analytical prove is provided, it seems that this equation has a unique solution with is approximately $(D_0,F_0)=(-0.56996,1.00781)$.

\section{Appendix: Computation of momenta}

The purpose of this section is to compute the momenta used to apply Theorem~\ref{thm:criticalitat-finite} on the two-parametric potential family and the Loud's family (see Proposition~\ref{prop:resultat-familia} and Section~\ref{sec:proof_ThmB}.) In both cases the momentum we are interested in is
\[
N\!\left[\bigl[\DD_{\nu}[f_{\mu}]\bigr]^1 \right]\ \text{ with }\nu=-1.
\]
Lemma~\ref{lem:lemes_varis} together with Definition~\ref{def:moment-finit} imply that
\[
N\!\left[\bigl[\DD_{\nu}[f_{\mu}]\bigr]^1 \right]=M\!\left[\bigl[ (\LL_{\nu}\circ\CC)[f_{\mu}]\bigr]_1\right].
\]
Using the formula in Lemma~\ref{lema:calcul-moment} on the right-hand side of the equality above, we get the equality
\[
N\!\left[\bigl[\DD_{\nu}[f_{\mu}]\bigr]^1 \right]=\lim_{R\rightarrow+\infty}\left(R^3\CC[f_{\mu}](R)-\int_0^R x^2\CC[f_{\mu}](x)dx \right).
\]
Finally, the change of variable $R=\phi^{-1}(r)=\frac{r}{\sqrt{1-r^2}}$ yields to
\begin{align*}
N\!\left[\bigl[\DD_{\nu}[f_{\mu}]\bigr]^1 \right]&=\lim_{r\rightarrow 1}\left(\frac{r^3f_{\mu}(r)}{\sqrt{1-r^2}}-\int_0^r \frac{z^2f_{\mu}(z)}{(1-z^2)^{\frac{3}{2}}}dz \right)\\
&=\lim_{r\rightarrow 1}\left(\frac{r^3f_{\mu}(r)}{\sqrt{1-r^2}}-\int_0^r \frac{f_{\mu}(z)}{(1-z^2)^{\frac{3}{2}}}dz+\int_0^r \frac{f_{\mu}(z)}{(1-z^2)^{\frac{1}{2}}}dz \right)\\
&=\lim_{r\rightarrow 1}\left(\frac{r^3f_{\mu}(r)}{\sqrt{1-r^2}}-\int_0^r \frac{f_{\mu}(z)}{(1-z^2)^{\frac{3}{2}}}dz\right)+N[f_{\mu}].
\end{align*}
Let $f_{\mu}(z)\!:=z\sqrt{h_0(\mu)}(g_{\mu}^{-1})''(z\sqrt{h_0(\mu)})-z\sqrt{h_0(\mu)}(g_{\mu}^{-1})''(-z\sqrt{h_0(\mu)})$ with $g_{\mu}\!:=\text{sgn}(x)\sqrt{V_{\mu}(x)}$. In both families, for parameters $\mu=\mu_0$ under consideration the momentum $N[f_{\mu_0}]$ vanishes. We point out that this fact is not a coincidence. Indeed, in both cases we are interested to bound the criticality of bifurcation parameters such that $\lim_{h\rightarrow h_0(\mu_0)}T_{\mu_0}'(h)=0.$ In other words, the bifurcation occurs in such a way that $\lim_{h\rightarrow h_0(\mu)}T_{\mu}'(h)$ is convergent for all $\mu\approx\mu_0$ and changes sign at $\mu=\mu_0$. From the definition of the period function, it turns out that
\begin{equation}\label{primera_constant}
N[f_{\mu}]=\lim_{h\rightarrow h_0(\mu)}\frac{T_{\mu}'(h)}{\sqrt{2}h_0(\mu)}
\end{equation}
and so we have $N[f_{\mu_0}]=0$ (we refer to~\cite[Corollary~3.12]{ManRojVil2016} for more details in this direction.) Hence, for parameters in such situation,
\begin{equation}\label{formula-momento2}
N\!\left[\bigl[\DD_{\nu}[f_{\mu}]\bigr]^1 \right]=\lim_{r\rightarrow 1}\left(\frac{r^3f_{\mu}(r)}{\sqrt{1-r^2}}-\int_0^r \frac{f_{\mu}(z)}{(1-z^2)^{\frac{3}{2}}}dz\right).
\end{equation}

Next two lemmas are useful for the forthcoming computations. See for instance $(6.2.1)$ and $(6.2.2)$ of~\cite{AS} for the first one. The first assertion of the second lemma was proved in~\cite[Lemma A.3]{ManRojVil2016}. The second assertion follows similarly.

\begin{lema}\label{lema:integral}
Let $\alpha$ and $\beta$ be any complex numbers with strictly positive real part. Then
\[
\int_0^1 u^{\alpha-1}(1-u)^{\beta-1}du = \int_0^{\infty} u^{\alpha-1}(1+u)^{-(\alpha+\beta)}du = \frac{\Gamma(\alpha)\Gamma(\beta)}{\Gamma(\alpha+\beta)},
\]
where $\Gamma$ denotes the Gamma function.
\end{lema}

\begin{lema}\label{lema:parts}
Let $\alpha$ and $\beta$ be real numbers such that $\alpha+\beta+1\neq 0$. Then
\[
\int u^{\alpha} (1+u)^{\beta} du = \frac{\beta}{\alpha+\beta+1}\int u^{\alpha}(1+u)^{\beta-1}du + \frac{1}{\alpha+\beta+1} u^{\alpha+1}(1+u)^{\beta},
\]
and
\[
\int u^{\alpha} (1-u)^{\beta} du = \frac{\alpha}{\alpha+\beta+1}\int u^{\alpha-1}(1-u)^{\beta}-\frac{1}{\alpha+\beta+1}u^{\alpha}(1-u)^{\beta+1}.
\]
\end{lema}

\begin{lema}\label{lema:integrals-phi}
Let $\Phi_{\mu}(x)\!:=h_0-V_{\mu}(x)$. The following hold:
\begin{enumerate}[$(a)$]
\item If $V_{\mu}$ with $\mu=(q,p)$ is defined by~\eqref{eq:potencial-familia} and $q\in(-1,-\tfrac{3}{5})$ then
\[
\int_0^s \Phi_{\mu}(x)^{-\frac{3}{2}}dx=\frac{2(p+1)^{\frac{3}{2}}(q+1)s^{-\frac{1}{2}(3q+1)}}{(p-q)(p+1)}u(s)^{-\frac{1}{2}} + \frac{2(p+1)^{\frac{3}{2}}(\lambda-1)}{(p-q)\bigl(\frac{p+1}{q+1}\bigr)^{\lambda}}\int_0^{k(s)}y^{\frac{1}{2}-\lambda}(1+y)^{\lambda-2}dy,
\]
and
\begin{align*}
\displaystyle\int_0^s \Phi_{\mu}(x)^{-\frac{5}{2}}dx &= \dfrac{2(p+1)^{\frac{3}{2}}(q+1)s^{-\frac{1}{2}(5q+3)}}{3(p-q)}u(s)^{-\frac{3}{2}}+\dfrac{4(p+1)^{\frac{1}{2}}(q+1)^2(\omega-1)s^{-\frac{1}{2}(5q+3)}}{3(p-q)}u(s)^{-\frac{1}{2}}\\
&\phantom{=}+\dfrac{4(p+1)^{\frac{5}{2}}(\omega-1)(\omega-2)}{3(p-q)\bigl(\frac{p+1}{q+1}\bigr)^{\omega}}\displaystyle\int_0^{k(s)} y^{\frac{3}{2}-\omega}(1+y)^{\mu-3}dy,
\end{align*}
where $u(s)\!:=\frac{p+1}{q+1}-s^{p-q}$, $k(s)\!:=\frac{1}{u(s)}s^{p-q}$, $\lambda\!:=\frac{1}{2}\frac{3p+1}{p-q}$ and $\omega\!:=\frac{1}{2}\frac{5p+3}{p-q}$.
\item If $V_{\mu}$ with $\mu=(D,F)$ is defined by~\eqref{V-loud} then
\[
\int_{-\frac{1}{F}}^s\!\Phi_{\mu}(x)^{-\frac{3}{2}}dx=\frac{(2-2F)^{\frac{3}{2}}(1-p_1)^{2F-\frac{3}{2}}}{D^{\frac{3}{2}}(p_2-p_1)^{\frac{3}{2}}}\left(2(v(s)^{-\frac{1}{2}}-1)+\!\int_{v(s)}^{1}\! y^{-\frac{3}{2}}((1-y)^{2-2F}(1-\alpha y)^{-\frac{3}{2}}-1)dy\right)
\]
and
\begin{align*}
\int_{-\frac{1}{F}}^s \Phi_{\mu}(x)^{-\frac{5}{2}}dx&=\frac{(2-2F)^{\frac{5}{2}}(1-p_1)^{4F-\frac{5}{2}}}{D^{\frac{5}{2}}(p_2-p_1)^{\frac{5}{2}}}\left( \frac{2}{3}(v(s)^{-\frac{3}{2}}-1)+(8F-8+5\alpha)(v(s)^{-\frac{1}{2}}-1)\right.\\
&\phantom{=}\left. + \int_{v(s)}^1 y^{-\frac{5}{2}}\left((1-y)^{4-4F}(1-\alpha y)^{-\frac{5}{2}}-1-(4F-4+\tfrac{5}{2}\alpha)y\right) dy \right)
\end{align*}
where $v(s)\!:=\frac{(Fs+1)^{-\frac{1}{F}}+p_1-1}{(Fs+1)^{-\frac{1}{F}}}$, $\alpha\!:=\frac{p_2-1}{p_2-p_1}$ and $p_1$, $p_2$ are defined in~\eqref{def:p1p2}.
\end{enumerate}
\end{lema}

\begin{prova}
Let us compute the integrals in $(a)$. Using the expression of $V_{\mu}$ in~\eqref{eq:potencial-familia} we have that
\[
\int_0^s \Phi(x)^{-\frac{3}{2}}dx=(p+1)^{\frac{3}{2}}\int_0^s x^{-\frac{3}{2}(q+1)}\left(\frac{p+1}{q+1}-x^{p-q}\right)^{-\frac{3}{2}}dx.
\]
We perform the change of variable $y=k(x)=\frac{x^{p-q}}{\frac{p+1}{q+1}-x^{p-q}}$. Then the integral above writes
\[
\int_0^s \Phi(x)^{-\frac{3}{2}}dx= \frac{(p+1)^{\frac{3}{2}}}{(p-q)\bigl(\frac{p+1}{q+1}\bigr)^{\lambda}}\int_0^{k(s)} y^{\frac{1}{2}-\lambda}(1+y)^{\lambda-1}dy,
\]
where $\lambda\!:=\frac{1}{2}\frac{3p+1}{p-q}$. This improper integral is divergent as $s$ tends to $\rho=\bigl(\frac{p+1}{q+1}\bigr)^{\frac{1}{p-q}}$, which will be the case under consideration. We shall make explicit the order of the singularity of this function using Lemma~\ref{lema:integral}. In order to apply it, we note that the power of $y$ in the integrand should be greater than $-1$ and the power of $(1+y)$ negative. The first condition is satisfied for $q\in(-1,-\tfrac{3}{5})$, $p>q$, but $\lambda-1=\frac{1+p+2q}{2(p-q)}$ is positive for some parameters. (In fact vanishes for the bifurcation parameters.) To overpass this situation, we apply Lemma~\ref{lema:parts}, getting
\[
\int_0^s \Phi(x)^{-\frac{3}{2}}dx=\frac{2(p+1)^{\frac{3}{2}}}{(p-q)\bigl(\frac{p+1}{q+1}\bigr)^{\lambda}}k(s)^{\frac{3}{2}-\lambda}(1+k(s))^{\lambda-1} + \frac{2(p+1)^{\frac{3}{2}}(\lambda-1)}{(p-q)\bigl(\frac{p+1}{q+1}\bigr)^{\lambda}}\int_0^{k(s)}y^{\frac{1}{2}-\lambda}(1+y)^{\lambda-2}dy.
\]
Using $u(s)\!:=\frac{p+1}{q+1}-s^{p-q}$ we can write $k(s)=s^{p-q}\frac{1}{u(s)}$ and $1+k(s)=\frac{p+1}{q+1}\frac{1}{u(s)}$ so,
\[
\int_0^s \Phi(x)^{-\frac{3}{2}}dx=\frac{2(p+1)^{\frac{3}{2}}(q+1)s^{-\frac{1}{2}(3q+1)}}{(p-q)(p+1)}u(s)^{-\frac{1}{2}} + \frac{2(p+1)^{\frac{3}{2}}(\lambda-1)}{(p-q)\bigl(\frac{p+1}{q+1}\bigr)^{\lambda}}\int_0^{k(s)}y^{\frac{1}{2}-\lambda}(1+y)^{\lambda-2}dy
\]
as desired. Similarly we can obtain the expression for the second integral in $(a)$. In this case, we need to apply Lemma~\ref{lema:parts} twice in order that the power of $(1+y)$ of the integrand is negative. We omit the details for the sake of brevity.

Let us now prove $(b)$. Using the expression of $V_{\mu}$ in~\eqref{V-loud} together with the change of variable $z=\phi^{-1}(x)=(Fx+1)^{-\frac{1}{F}}$ the first integral on item $(b)$ is written as
\[
\int_{\phi^{-1}(s)}^{+\infty} z^{2F-1} \bigl(h_0(\mu)-V_{\mu}(\phi(z))\bigr)^{-\frac{3}{2}}dz.
\]
By definition of $p_1$ and $p_2$ in~\eqref{def:p1p2} the difference on the integrand above decomposes as 
\[
h_0(\mu)-V_{\mu}(\phi(z))=\frac{D}{2-2F}(z-1+p_1)(z-1+p_2).
\]
Let us now consider the change of variable $y=v(x)=\frac{z+p_1-1}{z}=\frac{\phi^{-1}(x)+p_1-1}{\phi^{-1}(x)}.$ In this new variable, we have that
\[
\int_{-\frac{1}{F}}^s \Phi_{\mu}(x)^{-\frac{3}{2}}dx=\frac{(2-2F)^{\frac{3}{2}}(1-p_1)^{2F-\frac{3}{2}}}{D^{\frac{3}{2}}(p_2-p_1)^{\frac{3}{2}}}\int_{v(s)}^1 y^{-\frac{3}{2}}(1-y)^{2-2F}(1-\alpha y)^{-\frac{3}{2}}dy
\]
with $\alpha\!:=\frac{p_2-1}{p_2-p_1}$. As in $(a)$, we point out here that $v(s)$ tends to zero as $s$ tends to $x_r$ (see \eqref{ur_loud}) which will be the case under study. In particular, the above improper integral is divergent as $v(s)$ tends to zero. To deal with this situation, we add and subtract the value at $y=0$ of $(1-y)^{2-2F}(1-\alpha y)^{-\frac{3}{2}}$ inside the integral, which turns out to be exactly $1$. That is,
\[
\int_{-\frac{1}{F}}^s \Phi_{\mu}(x)^{-\frac{3}{2}}dx=\frac{(2-2F)^{\frac{3}{2}}(1-p_1)^{2F-\frac{3}{2}}}{D^{\frac{3}{2}}(p_2-p_1)^{\frac{3}{2}}}\left(2(v(s)^{-\frac{1}{2}}-1)+\int_{v(s)}^{1} y^{-\frac{3}{2}}((1-y)^{2-2F}(1-\alpha y)^{-\frac{3}{2}}-1)dy\right)
\]
as desired. Similarly we obtain the expression for the second integral in $(b)$. In this case, we need to add and subtract $1+(4F-4+\frac{5}{2}\alpha)y$ to the function $(1-y)^{2-2F}(1-\alpha y)^{-\frac{3}{2}}$, which turns out to be its Taylor's development up to degree one.
\end{prova}

\begin{lema}\label{lema:moment2}
Let $V_{\mu}$ be defined as in~\eqref{eq:potencial-familia} and $h_0(\mu)=\frac{p-q}{(p+1)(q+1)}$. If $\mu=(q,1)$ with $q\in(-1,-\tfrac{3}{5})$ and $\nu=-1$ then
\[
N\!\left[\bigl[\DD_{\nu}[f_{\mu}]\bigr]^1 \right]=\frac{2\sqrt{-2q}}{q(q+1)(3q+1)}\left(\frac{-2q}{q+1}\right)^{-\frac{1+5q}{1+3q}}\frac{\sqrt{\pi}\ \Gamma\left(-\frac{5q+3}{2(p-q)}\right)}{\Gamma\left(\frac{1}{2}-\frac{5q+3}{2(p-q)}\right)}.
\]
\end{lema}
\begin{prova}
For the parameters under consideration the equality deduced in~\eqref{formula-momento2} holds. For the sake of simplicity we shall omit the dependence on the parameter $\mu$ from now on. We stress that although the dependence on the parameters is not shown all the limits in the proof are uniform with respect to the parameters. Using the definition of $f$, that is
\[
f(z)=z\sqrt{h_0}(g^{-1})''(z\sqrt{h_0})-z\sqrt{h_0}(g^{-1})''(-z\sqrt{h_0})
\]
with $g(x)=\text{sgn}(x)\sqrt{V(x)}$, and the change of variable $z=g(x)/\sqrt{h_0}$, the following equality holds
\[
\int_0^r \frac{f(z)dz}{(1-z^2)^{\frac{3}{2}}}=2h_0\int_{g^{-1}(-r\sqrt{h_0})}^{g^{-1}(r\sqrt{h_0})} \frac{\eta_1(x)}{(h_0-V(x))^\frac{3}{2}}dx,
\]
with $\eta_1=\frac{1}{2}-\frac{VV''}{(V')^2}$. The function $g^{-1}(-r\sqrt{h_0})$ tends to zero as $r\rightarrow 1$. Moreover, $\eta_1(x)(h_0-V(x))^{-\frac{3}{2}}$ has a singularity at $x=0$ of order $x^{-\frac{5}{2}(q+1)}$, which is integrable since $q<-\frac{3}{5}$. Consequently,
\[
\lim_{r\rightarrow 1} \int_{g^{-1}(-r\sqrt{h_0})}^{g^{-1}(r\sqrt{h_0})} \frac{\eta_1(x)}{(h_0-V(x))^\frac{3}{2}}dx=\lim_{r\rightarrow 1} \int_{0}^{g^{-1}(r\sqrt{h_0})} \frac{\eta_1(x)}{(h_0-V(x))^\frac{3}{2}}dx.
\]
Let us define
\begin{equation}\label{funciones}
\Phi(x)\!:=h_0-V(x), \ l(x)\!:=-\frac{1}{V'(x)},\ h(x)\!:= h_0\Phi(x)^{-\frac{3}{2}}-\Phi(x)^{-\frac{1}{2}}.
\end{equation}
We have $\eta_1(x)(h_0-V(x))^{-\frac{3}{2}}=\frac{1}{2}\Phi(x)^{-\frac{3}{2}}-l'(x)h(x)$. Denoting $s=g^{-1}(r\sqrt{h_0})$ and integrating by parts, the previous discussion leads to the equality
\[
\int_0^r \frac{f(z)dz}{(1-z^2)^{\frac{3}{2}}}=h_0\int_0^s \Phi(x)^{-\frac{3}{2}}dx + 2h_0 \int_0^s l(x)h'(x)dx - 2h_0 l(s)h(s).
\]
Here we used that $\lim_{x\rightarrow 0} l(x)h(x)=0$, which is a direct computation. By definition, 
\[
l(x)h'(x)=-\frac{3}{2}h_0 \Phi(x)^{-\frac{5}{2}}+\frac{1}{2}\Phi(x)^{-\frac{3}{2}}.
\]
Then,
\begin{equation}\label{integral}
\int_0^r \frac{f(z)dz}{(1-z^2)^{\frac{3}{2}}}=2h_0\int_0^s \Phi(x)^{-\frac{3}{2}}dx-3h_0^2\int_0^s \Phi(x)^{-\frac{5}{2}}dx - 2h_0 l(s)h(s).
\end{equation}
From now on let us denote $u(s)\!:=\frac{p+1}{q+1}-s^{p-q}$. Using the expressions of $l$ and $h$ we have
\begin{equation}\label{lh}
l(s)h(s)=\frac{(p+1)^{\frac{1}{2}}s^{-\frac{1}{2}(3q+1)}}{s^{p-q}-1}u(s)^{-\frac{1}{2}}-\frac{(p-q)(p+1)^{\frac{1}{2}}s^{-\frac{1}{2}(5q+3)}}{(q+1)(s^{p-q}-1)}u(s)^{-\frac{3}{2}}.
\end{equation}
Substituting~\eqref{lh} and the integrals in assertion $(a)$ of Lemma~\ref{lema:integrals-phi} into~\eqref{integral}, and taking into account that
\[
\left(\frac{2h_0(p-q)(p+1)^{\frac{1}{2}}s^{-\frac{1}{2}(5q+3)}}{(q+1)(s^{p-q}-1)}-\frac{2h_0^2(p+1)^{\frac{3}{2}}(q+1)s^{-\frac{1}{2}(5q+3)}}{p-q}\right)u(s)^{-\frac{3}{2}}=\frac{2h_0(p+1)^{\frac{1}{2}}s^{-\frac{1}{2}(5q+3)}}{s^{p-q}-1}u(s)^{-\frac{1}{2}}
\]
we obtain
\begin{equation}\label{eq1}
\begin{array}{rl}
\displaystyle\int_0^r \frac{f(z)dz}{(1-z^2)^{\frac{3}{2}}}\!\!\!\! &=\dfrac{4h_0(p+1)^{\frac{3}{2}}(\lambda-1)}{(p-q)\bigl(\frac{p+1}{q+1}\bigr)^{\lambda}}\displaystyle\int_0^{k(s)}y^{\frac{1}{2}-\lambda}(1+y)^{\lambda-2}dy\\
&\phantom{=}-\dfrac{4h_0^2(p+1)^{\frac{5}{2}}(\omega-1)(\omega-2)}{(p-q)\bigl(\frac{p+1}{q+1}\bigr)^{\omega}}\displaystyle\int_0^{k(s)} y^{\frac{3}{2}-\omega}(1+y)^{\omega-3}dy\\
&\phantom{=}+2h_0(p+1)^{\frac{1}{2}}s^{-\frac{1}{2}(5q+3)}\left( \dfrac{2(q+1)s^{q+1}}{p-q}+\dfrac{1-s^{q+1}}{s^{p-q}-1}-\dfrac{2h_0(q+1)^2(\omega-1)}{p-q}\right)u(s)^{-\frac{1}{2}}.
\end{array}
\end{equation}
Let us now compute $\frac{r^3}{\sqrt{1-r^2}}f(r)$. Using again the expression of $f$ the previous function writes
\begin{equation}\label{eq-split}
\frac{r^3f(r)}{\sqrt{1-r^2}}=\frac{r^4\sqrt{h_0}(g^{-1})''(r\sqrt{h_0})}{\sqrt{1-r^2}}-\frac{r^4\sqrt{h_0}(g^{-1})''(-r\sqrt{h_0})}{\sqrt{1-r^2}}.
\end{equation}
First, the change of variable $s=g^{-1}(-r\sqrt{h_0})$ on the second quotient of the previous equality yields to
\[
\frac{r^4\sqrt{h_0}(g^{-1})''(-r\sqrt{h_0})}{\sqrt{1-r^2}}=\frac{2V(s)^2\bigl(V'(s)^2-2V(s)V''(s)\bigr)}{h_0V'(s)^3\sqrt{h_0-V(s)}}.
\]
The function $s=g^{-1}(-r\sqrt{h_0})\rightarrow 0$ as $r\rightarrow 1$ and the function in the right hand side of the previous equality has a singularity at $s=0$ of order $s^{-\frac{1}{2}(5q+3)}$. Since $q\in(-1,-\tfrac{3}{5})$ we have that
\[
\lim_{r\rightarrow 1} \frac{r^4\sqrt{h_0}(g^{-1})''(-r\sqrt{h_0})}{\sqrt{1-r^2}}=\lim_{s\rightarrow 0} \frac{2V(s)^2\bigl(V'(s)^2-2V(s)V''(s)\bigr)}{h_0V'(s)^3\sqrt{h_0-V(s)}}=0.
\] 
Second, the change of variable $s=g^{-1}(r\sqrt{h_0})$ yields to the same expression as before,
\[
\frac{r^4\sqrt{h_0}(g^{-1})''(r\sqrt{h_0})}{\sqrt{1-r^2}}=\frac{2V(s)^2\bigl(V'(s)^2-2V(s)V''(s)\bigr)}{h_0V'(s)^3\sqrt{h_0-V(s)}}.
\]
On account of the expression of $V$, expanding in series on $u=\frac{p+1}{q+1}-s^{p-q}$ we have that
\[
\frac{2V(s)^2\bigl(V'(s)^2-2V(s)V''(s)\bigr)}{h_0V'(s)^3\sqrt{h_0-V(s)}}=\frac{2h_0(p+1)^{\frac{1}{2}}s^{-\frac{1}{2}(5q+3)}}{(s^{p-q}-1)^3}\bigl( s^{q+1}(s^{p-q}-1)^2-2h_0(ps^{p-q}-q)\bigr)u^{-\frac{1}{2}} + o(u^{\frac{1}{2}}).
\]
Let us recall $\rho\!:=\left(\frac{p+1}{q+1}\right)^{\frac{1}{p-q}}$. At this point we claim that 
\begin{align*}
\lim_{s\rightarrow \rho}& \left( \frac{2h_0(p+1)^{\frac{1}{2}}s^{-\frac{1}{2}(5q+3)}}{(s^{p-q}-1)^3}\bigl( s^{q+1}(s^{p-q}-1)^2-2h_0(ps^{p-q}-q)\bigr)u(s)^{-\frac{1}{2}}\right.\\ &\ \ \left.-2h_0(p+1)^{\frac{1}{2}}s^{-\frac{1}{2}(5q+3)}\left( \dfrac{2(q+1)s^{q+1}}{p-q}+\dfrac{1-s^{q+1}}{s^{p-q}-1}-\dfrac{2h_0(q+1)^2(\omega-1)}{p-q}\right)u(s)^{-\frac{1}{2}}  \right)=0.
\end{align*}
Indeed, the previous expression writes
\[
2h_0(p+1)^{\frac{1}{2}}s^{-\frac{1}{2}(5q+3)}u^{-\frac{1}{2}}\left(\frac{2(q+1)s^{q+1}}{(p-q)(s^{p-q}-1)}u+\frac{2h_0(q+1)^2(\omega-1)}{p-q}-\frac{2h_0(ps^{p-q}-q)}{(s^{p-q}-1)^3}-\frac{1}{s^{p-q}-1}  \right).
\]
Using the expression of $h_0$ and $\omega$ appearing in Lemma~\ref{lema:integrals-phi},
\[
\frac{2h_0(q+1)^2(\omega-1)}{p-q}-\frac{2h_0(ps^{p-q}-q)}{(s^{p-q}-1)^3}-\frac{1}{s^{p-q}-1}= \frac{u}{(p-q)(p+1)}\left( c_1+c_2s^{p-q}+c_3 s^{2(p-q)}\right),
\]
for some $c_1,c_2,c_3\in\R[p,q]$. Substituting the previous equality on the limit above and taking into account that $u(s)=\frac{p+1}{q+1}-s^{p-q}\rightarrow 0$ as $s\rightarrow\rho$ the claim follows.
Therefore, 
\begin{align*}
N\!\left[\bigl[\DD_{\nu}[f_{\mu}]\bigr]^1 \right]&=\lim_{r\rightarrow 1}\left( \frac{r^3f(r)}{\sqrt{1-r^2}}-\int_0^r \frac{f(z)dz}{(1-z^2)^{\frac{3}{2}}}\right)\\
&= \lim_{s\rightarrow \rho} \left( \frac{2V(s)^2\bigl(V'(s)^2-2V(s)V''(s)\bigr)}{h_0V'(s)^3\sqrt{h_0-V(s)}}-\dfrac{4h_0(p+1)^{\frac{3}{2}}(\lambda-1)}{(p-q)\bigl(\frac{p+1}{q+1}\bigr)^{\lambda}}\displaystyle\int_0^{k(s)}y^{\frac{1}{2}-\lambda}(1+y)^{\lambda-2}dy\right. \\
&\phantom{= \lim_{s\rightarrow \rho}}+\dfrac{4h_0^2(p+1)^{\frac{5}{2}}(\omega-1)(\omega-2)}{(p-q)\bigl(\frac{p+1}{q+1}\bigr)^{\omega}}\displaystyle\int_0^{k(s)} y^{\frac{3}{2}-\omega}(1+y)^{\omega-3}dy\\
&\left.\phantom{= \lim_{s\rightarrow \rho}}-2h_0(p+1)^{\frac{1}{2}}s^{-\frac{1}{2}(5q+3)}\left( \dfrac{2(q+1)s^{q+1}}{p-q}+\dfrac{1-s^{q+1}}{s^{p-q}-1}-\dfrac{2h_0(q+1)^2(\omega-1)}{p-q}\right)u^{-\frac{1}{2}}\right)\\
&= \lim_{s\rightarrow \rho} \left( \dfrac{4h_0^2(p+1)^{\frac{5}{2}}(\omega-1)(\omega-2)}{(p-q)\bigl(\frac{p+1}{q+1}\bigr)^{\omega}}\displaystyle\int_0^{k(s)} y^{\frac{3}{2}-\omega}(1+y)^{\omega-3}dy\right.\\
&\phantom{= \lim_{s\rightarrow \rho}}\left.-\dfrac{4h_0(p+1)^{\frac{3}{2}}(\lambda-1)}{(p-q)\bigl(\frac{p+1}{q+1}\bigr)^{\lambda}}\displaystyle\int_0^{k(s)}y^{\frac{1}{2}-\lambda}(1+y)^{\lambda-2}dy\right).
\end{align*}
Let us now fix $p=-2q-1$. On account of the expression of $\lambda$ we have that $\lambda-1=0$. Finally, using the expressions of $\omega$ and $h_0$, the fact that $k(s)\rightarrow+\infty$ as $s\rightarrow+\infty$ and Lemma~\ref{lema:integral} we obtain
\[
N\!\left[\bigl[\DD_{\nu}[f_{\mu}]\bigr]^1 \right]=\frac{2\sqrt{-2q}}{q(q+1)(3q+1)}\left(\frac{-2q}{q+1}\right)^{-\frac{1+5q}{1+3q}}\frac{\sqrt{\pi}\ \Gamma\left(-\frac{5q+3}{2(p-q)}\right)}{\Gamma\left(\frac{1}{2}-\frac{5q+3}{2(p-q)}\right)}
\]
as we desired.
\end{prova}

\begin{lema}\label{lema:moment2-loud}
Let $V_{\mu}$ be defined as in~\eqref{V-loud}. The following holds:
\begin{enumerate}[$(a)$]
\item If $\mu=(D,F)$ satisfies $F\in(1,\tfrac{3}{2})$, $D\in(-\tfrac{3}{2},-\tfrac{1}{2})$ and $D+F>0$ then 
\[
N[f_{\mu}]=\frac{h_0(2-2F)^{\frac{3}{2}}(1-p_1)^{2F-\tfrac{3}{2}}}{D^{\frac{3}{2}}(p_2-p_1)^{\frac{3}{2}}}\frac{\sqrt{\pi}\Gamma(3-2F)}{\Gamma\bigl(\tfrac{5}{2}-2F\bigr)}{}_2F_1\bigl(\tfrac{3}{2},-\tfrac{1}{2};\tfrac{5}{2}-2F;\alpha).
\]
\item If $\mu=(D,F)$ satisfies $D=\mathcal{G}(F)$ and $F\in(1,\frac{5}{4})$ then 
\[
N\!\left[\bigl[\DD_{\nu}[f_{\mu}]\bigr]^1 \right]=-\frac{4h_0^2(2-2F)^{\frac{5}{2}}(1-p_1)^{4F-\frac{5}{2}}}{D^{\frac{5}{2}}(p_2-p_1)^{\frac{5}{2}}}\frac{\sqrt{\pi}\ \Gamma(5-4F)}{\Gamma\bigl(\tfrac{7}{2}-4F\bigr)}\ {}_2F_1\bigl(-\tfrac{3}{2},\tfrac{5}{2};\tfrac{7}{2}-4F;\alpha\bigr).
\]
\end{enumerate}
Here $\alpha\!:=\frac{p_2-1}{p_2-p_1}$, $h_0=\frac{F-D-1}{2F(F-1)(2F-1)}$ and $p_1$, $p_2$ are defined in~\eqref{def:p1p2}.
\end{lema}

\begin{prova}
The proof of the assertion in $(a)$ follows the same lines as the proof of $(b)$ using that
\[
N[f_{\mu}]=\lim_{r\rightarrow 1} \int_0^r \frac{f_{\mu}(z)dz}{(1-z^2)^{\frac{1}{2}}}.
\]
Since the proof of $(b)$ is richer in subtle technicalities, and for the sake of brevity, we decided to prove $(b)$ and omit the proof of assertion in $(a)$.

Let us show $(b)$. For the parameters $\mu=(D,F)$ with $D=\mathcal{G}(F)$ the equality in~\eqref{formula-momento2} holds. Again we omit the dependence on parameters for the sake of simplicity although all the limits are uniform with respect to the parameters. Following the same discussion as in the proof of Lemma~\ref{lema:moment2} we arrive to the equality
\begin{equation}\label{integral-loud}
\int_0^r \frac{f(z)dz}{(1-z^2)^{\frac{3}{2}}}=2h_0\int_{-\frac{1}{F}}^s \Phi(x)^{-\frac{3}{2}}dx-3h_0^2\int_{-\frac{1}{F}}^s \Phi(x)^{-\frac{5}{2}}dx - 2h_0 l(s)h(s),
\end{equation}
where the functions $\Phi$, $h$ and $l$ are defined in~\eqref{funciones}. The only difference lies on the interval of integration due to $g^{-1}(-r\sqrt{h_0})\rightarrow -\frac{1}{F}$ as $r\rightarrow 1$ in this case. 

Let us denote $u=u(s)\!:=\phi^{-1}(s)-1+p_1=(Fs+1)^{-\frac{1}{F}}-1+p_1.$ From the expressions of $h$ and $l$ in~\eqref{funciones} and the definition of $V$ and $V_1$ in~\eqref{V-loud} we have
\begin{equation}\label{lh-loud}
l(s)h(s)=\frac{(u+1-p_1)^{2F}(2-2F)^{\frac{1}{2}}}{D^{\frac{1}{2}}V_1(u+1-p_1)(u+p_2-p_1)^{\frac{1}{2}}} u^{-\frac{1}{2}} - \frac{h_0(u+1-p_1)^{4F}(2-2F)^{\frac{3}{2}}}{D^{\frac{3}{2}}V_1(u+1-p_1)(u+p_2-p_1)^{\frac{3}{2}}} u^{-\frac{3}{2}}.
\end{equation}
Moreover, equality~\eqref{eq-split} also holds and the limit of the second quotient of the right-hand side of the equality tends to zero as $r$ tends to one. Therefore
\[
\lim_{r\rightarrow 1} \frac{r^3f(r)}{\sqrt{1-r^2}}=\lim_{r\rightarrow 1}\frac{r^4\sqrt{h_0}(g^{-1})''(r\sqrt{h_0})}{\sqrt{1-r^2}}.
\]
The change of variable $s=g^{-1}(r\sqrt{h_0})$ yields to
\[
\frac{r^4\sqrt{h_0}(g^{-1})''(r\sqrt{h_0})}{\sqrt{1-r^2}}=\frac{2V(s)^2(V'(s)^2-2V(s)V''(s))}{h_0V'(s)^3\sqrt{h_0-V(s)}}.
\]
On account of the expressions in \eqref{V-loud}, if $z=\phi^{-1}(s)=(Fs+1)^{-\frac{1}{F}}$ we have
\[
\frac{r^4\sqrt{h_0}(g^{-1})''(r\sqrt{h_0})}{\sqrt{1-r^2}}=\frac{2z^{4F}(h_0-z^{-2F}V_0(z))^2(z^{-2F}V_1(z)^2-2(h_0-z^{-2F}V_0(z))V_2(z))}{h_0V_1(z)^3V_0(z)^{\frac{1}{2}}}.
\]
Finally, using the factorization $V_0(z)=\frac{D}{2-2F}(z-1+p_1)(z-1+p_2)$ and $u=u(s)=z-1+p_1$ defined above, after some algebraic manipulations we arrive to
\[
\frac{r^4\sqrt{h_0}(g^{-1})''(r\sqrt{h_0})}{\sqrt{1-r^2}}=\frac{2h_0(2-2F)^{\frac{1}{2}}\bigl((u+1-p_1)^{-2F}V_1(u+1-p_1)^2-2h_0V_2(u+1-p_1)\bigr)}{D^{\frac{1}{2}}V_1(u+1-p_1)^3(u+p_2-p_1)^{\frac{1}{2}}(u+1-p_1)^{-4F}} u^{-\frac{1}{2}}+o(u^{\frac{1}{2}}).
\]
At this point we substitute the previous equality and~\eqref{integral-loud} into~\eqref{formula-momento2}, and we use the equality~\eqref{lh-loud} and the integrals in item $(b)$ of Lemma~\ref{lema:integrals-phi} to have an explicit expression for the momentum under consideration. On account that $v(s)=\frac{\phi^{-1}(s)-1+p_1}{\phi^{-1}(s)}=\frac{u}{u+1-p_1}$ (see Lemma~\ref{lema:integrals-phi}), we can collect the expression of the momentum as follows
\begin{align*}
N\!\left[\bigl[\DD_{\nu}[f_{\mu}]\bigr]^1 \right]&=\lim_{r\rightarrow 1}\left( \frac{r^3f(r)}{\sqrt{1-r^2}}-\int_0^r \frac{f(z)dz}{(1-z^2)^{\frac{3}{2}}}\right)\\
&=\lim_{u\rightarrow 0}\left(a(u)u^{-\frac{3}{2}}+b(u)u^{-\frac{1}{2}}+\frac{2h_0(2-2F)^{\frac{3}{2}}(1-p_1)^{2F-\frac{3}{2}}}{D^{\frac{3}{2}}(p_2-p_1)^{\frac{3}{2}}} I_1(u)\right. \\
&\phantom{=\lim_{u\rightarrow 0}}\ \ \left.+\frac{3h_0^2(2-2F)^{\frac{5}{2}}(1-p_1)^{4F-\frac{5}{2}}}{D^{\frac{5}{2}}(p_2-p_1)^{\frac{5}{2}}}I_2(u)\right),
\end{align*}
where
\begin{align*}
I_1(u)\!&:=2-\int_{\frac{u}{u+1-p_1}}^1 y^{-\frac{3}{2}}((1-y)^{2-2F}(1-\alpha y)^{-\frac{3}{2}}-1)dy,\\
I_2(u)\!&:=-\frac{22}{3}+8F+5\alpha- \int_{\frac{u}{u+1-p_1}}^1 y^{-\frac{5}{2}}\left((1-y)^{4-4F}(1-\alpha y)^{-\frac{5}{2}}-1-(4F-4+\tfrac{5}{2}\alpha)y\right) dy,
\end{align*}
and both $a$ and $b$ are analytic functions at $u=0$. Here $\alpha=\frac{p_2-1}{p_2-p_1}$. First, from the first part of the proof of Lemma~\ref{prop:resultat-familia} we know that $N\!\left[\bigl[\DD_{\nu}[f_{\mu}]\bigr]^1 \right]$ exists and so the previous limit is finite. (This can also be deduced from the fact that at the moment when certain momentum in Theorem~\ref{thm:criticalitat-finite} needs to be computed, such momentum exists and it is finite.) The previous fact, together with the analyticity of $a$ and $b$, implies that $a(u)u^{-\frac{3}{2}}+b(u)u^{-\frac{1}{2}}\rightarrow 0$ as $u\rightarrow 0$. Second, we point out that $I_1(0)$ is the expression of the curve $D=\mathcal{G}(F)$ found in~\cite{MarMarVil2006} which vanishes for bifurcation parameters we are considering. Consequently, taking the limit on the equality above,
\begin{align*}
N\!\left[\bigl[\DD_{\nu}[f_{\mu}]\bigr]^1 \right]=\frac{3h_0^2(2-2F)^{\frac{5}{2}}(1-p_1)^{4F-\frac{5}{2}}}{D^{\frac{5}{2}}(p_2-p_1)^{\frac{5}{2}}}I_2(0).
\end{align*}
The result will follow once we prove that
\[
I_2(0)=-\frac{4\sqrt{\pi}\ \Gamma(5-4F)}{3\ \Gamma\bigl(\tfrac{7}{2}-4F\bigr)}\ {}_2F_1\bigl(-\tfrac{3}{2},\tfrac{5}{2};\tfrac{7}{2}-4F;\alpha\bigr).
\]
To do so we shall show that the power series of the functions at both sides of the equality coincides. We notice that for parameters under consideration, $\abs{\alpha}<1$ and so the Hypergeometric function is analytic and well defined as a power series (see~\cite[Section 15]{AS}). First notice that using the Newton's binomial the integrand of $I_2(0)$ writes
\begin{equation}\label{newton}
\frac{-1+(1-y)^{4-4F}+(4-4F)y}{y^{\frac{5}{2}}}+
\frac{5(1-y)^{4-4F}-5}{2y^{\frac{3}{2}}}\alpha+(1-y)^{4-4F}y^{-\frac{5}{2}}\sum_{k=2}^{\infty}\left(\begin{matrix}
\frac{3}{2}+k\\ \frac{3}{2}
\end{matrix}\right)\alpha^k y^k.
\end{equation}
By the assertion in Lemma~\ref{lema:integral} and on account that $F<5/4$, for each $k\geq 2$ we have 
\[
\int_0^1 (1-y)^{4-4F}y^{k-\frac{5}{2}}dy=\frac{\Gamma(5-4F)\Gamma\bigl(-\tfrac{3}{2}+k\bigr)}{\Gamma\bigl(\tfrac{7}{2}-4F+k\bigr)}.
\]
Since the function inside the integral is positive for each $k\geq 2$,  Tonelli's theorem states that summation and integral signs interchange. Thus we obtain
\begin{equation}\label{hyper1}
\begin{array}{rl}
\displaystyle\int_0^1 \sum_{k=2}^{\infty}\left(\begin{matrix}
\frac{3}{2}+k\\ \frac{3}{2}
\end{matrix}\right)\alpha^k(1-y)^{4-4F}y^{k-\frac{5}{2}}dy\!\!\!&=\displaystyle\sum_{k=2}^{\infty}\left(\begin{matrix}
\frac{3}{2}+k\\ \frac{3}{2}
\end{matrix}\right)\frac{\Gamma(5-4F)\Gamma\bigl(-\tfrac{3}{2}+k\bigr)}{\Gamma\bigl(\tfrac{7}{2}-4F+k\bigr)}\alpha^k\\
&=\displaystyle\sum_{k=2}^{\infty}\frac{\Gamma(5-4F)\Gamma\bigl(-\tfrac{3}{2}+k\bigr)\Gamma\bigl(\tfrac{5}{2}+k\bigr)}{\Gamma\bigl(\tfrac{7}{2}-4F+k\bigr)\Gamma\bigl(\tfrac{5}{2}\bigr)}\frac{\alpha^k}{k!}.
\end{array}
\end{equation}
Let us now integrate the second term in~\eqref{newton}. The integration of the first term follows similarly and we omit the computation for the sake of shortness. Taking $\epsilon>0$ small we have
\[
\int_{\epsilon}^1 \frac{5(1-y)^{4-4F}-5}{2y^{\frac{3}{2}}}\alpha dy= \frac{5}{2}\alpha\int_{\epsilon}^1 y^{-\frac{3}{2}}(1-y)^{4-4F}dy+5\alpha (1-\epsilon^{-\frac{1}{2}}).
\]
With the aim of applying Lemma~\ref{lema:integral} in view we need the powers of $y$ and $(1-y)$ to be greater than $-1$. Indeed, since $F\in(1,\tfrac{5}{4})$ we already have $4-4F>-1$. On the contrary, the power of $y$ do not satisfy the assumptions of the lemma. To overcome this we use Lemma~\ref{lema:parts} and we get
\[
\int_{\epsilon}^1 y^{-\frac{3}{2}}(1-y)^{4-4F}dy=-2\bigl(\tfrac{9}{2}-4F\bigr)\int_{\epsilon}^1 y^{-\frac{1}{2}}(1-y)^{4-4F}dy+2\epsilon^{-\frac{1}{2}}(1-\epsilon)^{5-4F}.
\]
Substituting this equality on the previous expression, tending $\epsilon$ to zero and using Lemma~\ref{lema:integral} we get
\[
\int_{0}^1 \frac{5(1-y)^{4-4F}-5}{2y^{\frac{3}{2}}}\alpha dy= 5\alpha-5\alpha\bigl(\tfrac{9}{2}-4F\bigr)\int_{0}^1 y^{-\frac{1}{2}}(1-y)^{4-4F}dy=5\alpha-\frac{5\alpha\sqrt{\pi}\Gamma(5-4F)}{\Gamma\bigl(\tfrac{9}{2}-4F\bigr)}.
\]
As we noticed before, the same procedure shows that
\[
\int_0^1 \frac{-1+(1-y)^{4-4F}+(4-4F)y}{y^{\frac{5}{2}}}dy = -\frac{22}{3}+8F+\frac{4\sqrt{\pi}\Gamma(5-4F)}{3\Gamma\bigl(\tfrac{7}{2}-4F\bigr)}.
\]
In this last case Lemma~\ref{lema:parts} must be applied twice. Summing these two last expressions together with~\eqref{hyper1} we have that
\begin{align*}
I_2(0)&=-\frac{4\sqrt{\pi}\Gamma(5-4F)}{3\Gamma\bigl(\tfrac{7}{2}-4F\bigr)}+\frac{5\alpha\sqrt{\pi}\Gamma(5-4F)}{\Gamma\bigl(\tfrac{9}{2}-4F\bigr)}-\sum_{k=2}^{\infty}\frac{\Gamma(5-4F)\Gamma\bigl(-\tfrac{3}{2}+k\bigr)\Gamma\bigl(\tfrac{5}{2}+k\bigr)}{\Gamma\bigl(\tfrac{7}{2}-4F+k\bigr)\Gamma\bigl(\tfrac{5}{2}\bigr)}\frac{\alpha^k}{k!}\\
&=-\frac{4\sqrt{\pi}\Gamma(5-4F)}{3\Gamma\bigl(\tfrac{7}{2}-4F\bigr)}\left(1+\frac{15\alpha}{2(8F-7)}+\sum_{k=2}^{\infty}\frac{\Gamma\bigl(-\tfrac{3}{2}+k\bigr)\Gamma\bigl(\tfrac{5}{2}+k\bigr)\Gamma\bigl(\tfrac{7}{2}-4F\bigr)}{\Gamma\bigl(-\tfrac{3}{2}\bigr)\Gamma\bigl(\tfrac{5}{2}\bigr)\Gamma\bigl(\tfrac{7}{2}-4F+k\bigr)}\frac{\alpha^k}{k!}\right)\\
&=-\frac{4\sqrt{\pi}\Gamma(5-4F)}{3\Gamma\bigl(\tfrac{7}{2}-4F\bigr)}\ {}_2F_1\left(-\tfrac{3}{2},\tfrac{5}{2};\tfrac{7}{2}-4F;\alpha\right),
\end{align*}
where on the first equality we use common properties of the Gamma function and on the second equality we use the definition of the Hypergeometric function ${}_2F_1$ as a power series on account that $\abs{\alpha}<1$.
\end{prova}

\begin{rem}
The Hypergeometric function ${}_2F_1\bigl(a,b;c;z)$ can be continued analytically for any complex number with $\abs{z}\geq 1$ along any path in the complex plane that avoids the branch points 1 and infinity. For parameters $\mu=(D,F)$ satisfying $F\in(1,\tfrac{3}{2})$, $D\in(-\tfrac{3}{2},-\tfrac{1}{2})$ and $D+F>0$ the property $\abs{\alpha(\mu)}<1$ do not hold. However, it holds that $\alpha(\mu)<1$. Thus the function
\[
\frac{1}{\Gamma\bigl(\tfrac{5}{2}-2F\bigr)}{}_2F_1\bigl(\tfrac{3}{2},-\tfrac{1}{2};\tfrac{5}{2}-2F;\alpha(\mu))
\]
in the statement of Lemma~\ref{lema:moment2-loud} is well defined as a power series.
\end{rem}

\begin{prooftext}{Proof of Proposition~\ref{prop:corba-bif-expr}}

As it is shown in~\eqref{primera_constant} the bifurcation curve $D=\mathcal{G}(F)$ must coincide with the zero set of $N[f_{\mu}]$. For the parameters under consideration, the first assertion in Lemma~\ref{lema:moment2-loud} shows that this zero set coincide with the zero set of the function
${}_2F_1\bigl(\tfrac{3}{2},-\tfrac{1}{2};\tfrac{5}{2}-2F;\alpha)$ where $\alpha\!:=\frac{p_2-1}{p_2-p_1}$. This proves the first assertion of the result. To show the second assertion notice that 
\[
\lim_{(D,F)\rightarrow (-1,5/4)} {}_2F_1\bigl(\tfrac{3}{2},-\tfrac{1}{2};\tfrac{5}{2}-2F;\alpha(\mu))=1
\]
since $\alpha(\mu)\rightarrow 0$ as $D\rightarrow-1$. Therefore the result follows by the limit
$\lim_{(D,F)\rightarrow(-1,5/4)} \frac{\Gamma(3-2F)}{\Gamma\bigl(\tfrac{5}{2}-2F\bigr)}=0.$
\end{prooftext}

\section*{Acknowledgements}
The author thanks D. Marin and J. Villadelprat for fruitful discussions and valuable comments during the development of this work. The author is partially supported by the MINECO/FEDER grant MTM2017-82348-C2-1-P and the  MINECO/FEDER grant MTM2017-86795-C3-1-P.


\end{document}